\documentclass{amsart}
\usepackage{amsfonts,amssymb,amsmath,amsthm}
\usepackage{url}
\usepackage{enumerate}

\theoremstyle{definition}

\numberwithin{equation}{section}
\email{karkri.rafik@gmail.com}
\email{hzoubeir2014@gmail.com}
\input{tcilatex}

\begin{document}
\title[About Schauder frames and besselian Schauder frames]{About Schauder
frames and besselian Schauder frames of Banach spaces.}
\author[R. Karkri]{Rafik Karkri}
\address{Ibn Tofail University, Departement of Mathematics, Faculty of
sciences, P.B.133, Kenitra, Morocco. }
\author[H. Zoubeir ]{Hicham Zoubeir}
\address{Ibn Tofail University, Departement of Mathematics, Faculty of
sciences, P.B.133, Kenitra, Morocco. }

\begin{abstract}
In this paper, we introduce, for a separable Banach spacea new notion of
besselian paires and of besselian Schauder frames for which we prove for
some fundamental results some of which will allow us to obtain :

(i) a theorem providing necessary and sufficient conditions for the
existence of a Schauder frame for a separable Banach space;

(ii) a theorem providing necessary conditions for the existence of a
besselian Schauder frame for a weakly sequentially complete separable Banach
space.

From the conditions (ii) we deduce that the Banach space $L_{1}\left(
[0,1]\right) $ don't have a besselian Schauder frame. Finally we will show,
the existence of two universal Banach spaces in terms of which conditions
(i) and (ii) can be expressed.
\end{abstract}

\maketitle

\section{Introduction}

In $1946,$ Gabor \cite{gab} performed a new method for the signal
reconstruction from elementary signals. In $1952,$ Duffin and Schaeffer \cite%
{duf} developped, in the field of nonharmonic series, a similar tool and
introduced frame theory for Hilbert spaces. For more than thirty years, the
results of Duffin and Schaeffer has not received from the mathematical
community, the interest they deserve, until the publication of the work of
Young \cite{you} where \ the author studied frames in abstract Hilbert
spaces. In $1986,$ the work of Daubechies, Grossmann and Meyer gave to frame
\ theory the momentum it lacked and allowed it to be widely studied. This
contributed, among other things, to the wider developpement\ of wavelet
theory. The concept of atomic decompositions was introduced, in $1988,$ by
Feichtinger and Gr\"{o}chenig \cite{fei}, in order to extend the definition
of frames from the seting of Hilbert spaces to that of general separable
Banach spaces. In $1991,$ Gr\"{o}chenig \cite{gro}, presented a
generalisation of the notions of atomic decomposition and of synthesis
operator and introduced the definition of Banach frames. In $2001,$
Aldroubi, Sun and Tang \cite{ald} introduced the concepts of $p$-frames. In $%
2003,$ Christensen and Stoeva \cite{chr1} extended the definition of $p$%
-frames, by replacing the sequence space $l^{p}$ by a more general scalar
sequence space $X_{d}.$ By "getting rid" of the sequence spaces $X_{d}$ in
the definition of atomic decompositions, Cassaza, Han, Larson in $1999$ \cite%
{lar1} and in $2000$ Han and Larson \cite{lar2}, generalized the notion of
atomic decompositions by introducing the new notion of Schauder frames. One
of the peculiarities of Schauder frames is that they constitute a natural
extension of the concept of Schauder basis.

On the other hand, we derive from the definitions of frames and $p$-frames
respectively, a new classes of sequences, the class Bessel sequences and the
class of $p$-Bessel sequences respectively, for which there exists an
extensive literature with various lines of research :\  \  \  \  \  \  \  \  \  \  \  \
\  \  \  \  \  \  \  \  \  \  \  \  \  \  \  \  \  \  \  \  \  \  \  \  \  \  \  \  \  \  \  \  \  \  \  \  \  \
\  \ 

In $2006,$ Jia \cite{jia} studied Bessel sequences in Sobolev spaces. In $%
2008,$ Casazza and Leonhard \cite{cas3} proved that any Bessel sequence in a
finite-dimensional space can be extended to a tight frame. This result was
extended, in $2009,$ by Li and Sun \cite{li}, to finite-dimensional spaces.
In the same year Christensen, Kim and Kim \cite{chr2} proved that it is
possible, in any separable Hilbert space to extend each pair of Bessel
sequences to a pair of mutually dual frames. In $2010,$ Rahimi and Balazs 
\cite{rah} investigated multipliers for $p$-Bessel sequences in Banach
spaces. In $2015,$ Bakic and Beric \cite{bak} studied finite extensions of
Bessel seaquences in infinite dimensional Hilbert spaces. In $2015,$ Arias,
Corach and Pacheco \cite{ari} studied the set of all Bessel sequences for a
separable Hilbert as a Banach space. In $2015,$ Dehgan and Mesbah \cite{deh}
endowed the set of all Bessel sequences for a separable Hilbert with a
structure of $C^{\ast }$ algebra.

In this paper, we introduce a new notion of besselian paires and of
besselian Schauder frames of a Banach space.

Let us explain briefly how we "discovered" theses notions :

\begin{itemize}
\item The starting point of our reasoning was the following observation
concerning Bessel sequences of separable Hilbert spaces : we can prove (see
lemma 4.2. in section 4., with $p=2$) that a sequence $\left( x_{n}\right)
_{n\in 
\mathbb{N}
^{\ast }}$ of a separable Hilbert space $\left( \mathcal{H},\left \langle
\cdot ,\cdot \right \rangle \right) $ is a Bessel sequence if and only if
there exists a constant $C>0$ such that :%
\begin{equation}
\underset{n=1}{\overset{+\infty }{\sum }}\left \vert x_{n}^{\ast }(x)\right
\vert \left \vert y^{\ast }\left( x_{n}\right) \right \vert \leq C\left
\Vert x\right \Vert _{\mathcal{H}}\left \Vert y^{\ast }\right \Vert _{%
\mathcal{H}^{\ast }},\text{ }x\in \mathcal{H},\text{ }y^{\ast }\in \mathcal{H%
}^{\ast }  \label{one}
\end{equation}%
where $\mathcal{H}^{\ast }$ denotes the dual space of $\mathcal{H},$ and $%
x_{n}^{\ast }$ is for every $n\in 
\mathbb{N}
^{\ast }$the element of $\mathcal{H}^{\ast }$ defined$,$ for each $x\in 
\mathcal{H}$, by the relation $x_{n}^{\ast }(x)=\left \langle
x,x_{n}\right
\rangle .$

\item The following steps in our approach was :
\end{itemize}

\begin{enumerate}
\item To replace the separable Hilbert space by a general separable banach
space $X.$

\item To "get rid" of the relation, in the assertion (\ref{one}), between
the vector $x_{n}$ and the linear form $x_{n}^{\ast }$ and to consider
instead of a sequence of vectors $\left( x_{n}\right) _{n\in 
\mathbb{N}
^{\ast }}$ a sequence $\left( \left( x_{n},y_{n}^{\ast }\right) \right)
_{n\in 
\mathbb{N}
^{\ast }}$ of elements of $X\times X^{\ast }$ that we called a paire of $X.$
Hence considering a paire $\mathcal{F}:=\left( \left( x_{n},y_{n}^{\ast
}\right) \right) _{n\in 
\mathbb{N}
^{\ast }}$ of a Banach space $X,$ we say that $\mathcal{F}$ is a besselian
paire of $X$ if the following condition holds :%
\begin{equation}
\underset{n=1}{\overset{+\infty }{\sum }}\left \vert y_{n}^{\ast }(x)\right
\vert \left \vert y^{\ast }\left( x_{n}\right) \right \vert \leq C\left
\Vert x\right \Vert _{\mathcal{H}}\left \Vert y^{\ast }\right \Vert _{%
\mathcal{H}^{\ast }},\text{ }x\in E,\text{ }y^{\ast }\in E^{\ast }
\label{two}
\end{equation}%
where $C>0$ is a constant.

\item Since Schauder frames on the Banach space $X$ are also paires of $X$,
we say that a Schauder frame $\mathcal{F}$ is a besselian Schauder frame of $%
X$ if it is a besselian paire in the sens of definition \ref{two}.
\end{enumerate}

The paper is structured as follows. After stating, in section 2, the main
necessary definitions and notations, we devote section 3 to some examples of
classical separable banach spaces (mostly scalar sequence spaces) with
besselian Schauder frames. In section 4, we prove for besselian Schauder
frame some fundamental results, some of which will allow us to obtain in
section 5 :

(i) a theorem providing necessary and sufficient conditions for the
existence of a Schauder frame for a separable Banach space;

(ii) a theorem providing necessary conditions for the existence of a
besselian Schauder frame for a weakly sequentially complete separable Banach
space. As a corollary\ of the second theorem, we show that the Banach space $%
L_{1}\left( [0,1]\right) $ don't have a besselian Schauder frame. Finally we
will show, in section 5, the existence of two universal separable Banach
spaces in terms of which conditions (i) and (ii) can be expressed.

Let us point out that we were inspired, in the proof of the theorems of
section 5, by the work \cite{pel3} of A. Pe\l czy\'{n}ski and that the proof
of lemma 4.2. is a generalisation of the proof of a part of proposition 4 in
(\cite{woj}, page 59).

\section{Main definitions and notations}

Let $X$ be a sparable Banach space on $\mathbb{K}\in \left \{ \mathbb{R},%
\mathbb{C}\right \} $, $X^{\ast }$ its dual, $\mathbb{I}$ the set $\mathbb{N}%
^{\ast }$or a set of the form $\left \{ 1,...,r\right \} $ where $r\in 
\mathbb{N}^{\ast }$and $\left( \left( x_{n},y_{n}^{\ast }\right) \right)
_{n\in \mathbb{I}}$ a sequence of elements of $X\times X^{\ast }.$

\begin{enumerate}
\item We make the \ convention that each sum indexed by an expression of the
form $j\in \emptyset ,$ is equal to zero.

\item For each $q\in ]1,+\infty \lbrack $, we set $q^{\ast }=\dfrac{q}{q-1}.$

\item For each $\left( j,k\right) \in 
\mathbb{N}
^{2},$ we denote by $\delta _{jk}$ the so-called Kroneker symbole defined by 
\begin{equation*}
\delta _{jk}:=\left \{ 
\begin{array}{c}
1\text{ if }j=k \\ 
0\text{ else}%
\end{array}%
\right.
\end{equation*}

\item We denote by $\mathcal{P}_{\mathbb{N}^{\ast }}$ (resp. $\mathcal{D}_{%
\mathbb{N}^{\ast }}$) the set of all the nonempty subsets (resp. all the
nonempty finite subsets) of $\mathbb{N}^{\ast }.$ For each $G\in \mathcal{P}%
_{\mathbb{N}^{\ast }}$ we denote by $\mathcal{P}_{G}$ the set of all the
nonempty subsets of $G.$

\item For each $A\in \mathcal{P}_{\mathbb{N}^{\ast }},$ we denote by $\min
\left( A\right) $ the minimal element of $A.$We denote by $\mathfrak{S}_{%
\mathbb{N}^{\ast }}$ the set of all the permutations of $\mathbb{N}^{\ast }.$

\item For each mapping $F:G\rightarrow H$ between nonempty sets, and for
each nonempty subset $L$ of $G,$ we denote by $F_{\left \vert L\right. }$
the restriction of $F$ to $L.$

\item We denote by $0_{X}$ the null vector of $X.$

\item Let $\left( x_{n}\right) _{n\in \mathbb{I}}$ be a sequence of elements
of $E.$ If $\mathbb{I}$ is of the form $\mathbb{I}=\left \{ 1,...,r\right \} 
$ $\left( r\in 
\mathbb{N}
^{\ast }\right) $ then the notation $\sum_{n\in \mathbb{I}}x_{n}$ will
represent, of course, the sum $\underset{n=1}{\overset{r}{\sum }}x_{n}.$ If $%
\mathbb{I}=%
\mathbb{N}
^{\ast }$ and the series $\sum x_{n}$ is convergent in $X$, then the
notation $\sum_{n\in \mathbb{I}}x_{n}$ will represent, of course, the sum $%
\underset{n=1}{\overset{+\infty }{\sum }}x_{n}$ of the series $\sum x_{n}.$

\item We denote by $\mathbb{B}_{X\text{ \ }}$the closed unit ball of $X:$%
\begin{equation*}
\mathbb{B}_{X}:=\left \{ x\in X:\left \Vert x\right \Vert _{X}\leq 1\right \}
\end{equation*}

\item We denote by $L\left( X\right) $ the set of all\ bounded linear
operators $f:X\rightarrow X$ . It is well-known that $L\left( X\right) $ is
a Banach space for the norm $\left \Vert \cdot \right \Vert _{L\left(
X\right) }$ defined by the formula :%
\begin{equation*}
\left \Vert f\right \Vert _{L\left( X\right) }:=\underset{x\in \mathbb{B}_{X}%
}{\sup }\left \Vert f\left( x\right) \right \Vert _{X}
\end{equation*}

\item Let $\left( x_{n}\right) _{n\in 
\mathbb{N}
^{\ast }}$ be a sequence of elements of $X.$ The series $\sum x_{n}$ of s of 
$X$ is said to be weakly unconditionally convergent if for every functional $%
x^{\ast }\in X^{\ast }$ the scalar series $\sum \left \vert x^{\ast }\left(
x_{n}\right) \right \vert $ is convergent.

\item The Banach space $X$ is called weakly sequentially complete if for
each sequence $(x_{n})_{n\in \mathbb{N}^{\ast }}$ of $X$ such that $\underset%
{n\rightarrow +\infty }{\lim }y^{\ast }(x_{n})$ exists for every $y^{\ast
}\in X^{\ast },$ there exists $x\in X$ such that $\underset{n\rightarrow
+\infty }{\lim }y^{\ast }(x_{n})=y^{\ast }(x)$ for every $y^{\ast }\in
X^{\ast }$.

\item Let $\left( \mathcal{H},\left \langle .,.\right \rangle \right) $ be a
Hilbert space on $\mathbb{K}$ and $f$ an element of $\mathcal{H}$. We denote
by $T_{f}$ the continuous linear form $T_{f}:\mathcal{H}\longrightarrow 
\mathbb{K}$ defined by: $T_{f}(x)=\left \langle x,f\right \rangle ,\; \;x\in 
\mathcal{H}$.

\item Let $\left( \mathcal{H},\left \langle .,.\right \rangle \right) $ be a
Hilbert space and $\left( y_{n}\right) _{n\in \mathbb{N}^{\ast }}$ be a
sequence of elements of $\mathcal{H}$. $\left( y_{n}\right) _{n\in \mathbb{N}%
^{\ast }}$ is called a Bessel sequence if there exists a constant $B>0$ such
that: 
\begin{equation*}
\underset{n=1}{\overset{+\infty }{\sum }}\left \vert \left \langle
y,y_{n}\right \rangle \right \vert ^{2}\leq B\left \Vert y\right \Vert _{%
\mathcal{H}}^{2},\; \;y\in \mathcal{H}
\end{equation*}

\item Let $\left( x_{n}\right) _{n\in 
\mathbb{N}
^{\ast }}$ be a sequence of elements of $X$ and $q\in \left] 1,+\infty %
\right[ $. We say that the sequence $\left( x_{n}\right) _{n\in 
\mathbb{N}
^{\ast }}$ is weakly $q$-summable in $X$ if the series $\sum \left \vert
f^{\ast }\left( x_{n}\right) \right \vert ^{q}$ is convergent for each $%
f^{\ast }\in X^{\ast }.$

\item Let $\left( f_{n}^{\ast }\right) _{n\in 
\mathbb{N}
^{\ast }}$ be a sequence of elements of $X^{\ast }$. We say that the
sequence $\left( f_{n}^{\ast }\right) _{n\in 
\mathbb{N}
^{\ast }}$ is *-weakly $q$-summable in $X^{\ast }$ if the series $\sum
\left
\vert f_{n}^{\ast }\left( x\right) \right \vert ^{q}$ is convergent
for each $x\in X$.

\item The sequence $\left( \left( x_{n},y_{n}^{\ast }\right) \right) _{n\in 
\mathbb{I}}$ is called a paire of $X$.

\item A paire of $X$ of the form $\left( \left( x_{n},y_{n}^{\ast }\right)
\right) _{n\in \left \{ 1,...,r\right \} }$ is called a Schauder frame of $X$
if for all $x\in X$, the following formula holds :%
\begin{equation*}
\underset{n=1}{\overset{r}{\sum }}y_{n}^{\ast }\left( x\right) x_{n}=x
\end{equation*}

\item A paire of $X$ of the form $\left( \left( x_{n},y_{n}^{\ast }\right)
\right) _{n\in \mathbb{N}^{\ast }}$ is called a Schauder frame (resp.
unconditional Schauder frame ) of $X$ if for all $x\in X$, the serie $\sum
y_{n}^{\ast }\left( x\right) x_{n}$ is convergent (resp. unconditionally
convergent) in $X$ to $x$.

\item The paire $\left( \left( x_{n},y_{n}^{\ast }\right) \right) _{n\in 
\mathbb{I}}$ of $X$ is said to be a besselian paire of $X$ if there exists a
constant $A>0$ such that 
\begin{equation*}
\underset{n\in \mathbb{I}}{\sum }\left \vert y_{n}^{\ast }\left( x\right)
\right \vert \left \vert y^{\ast }\left( x_{n}\right) \right \vert \leq
A\left \Vert x\right \Vert _{X}\left \Vert y^{\ast }\right \Vert _{X^{\ast }}
\end{equation*}%
for each $\left( x,y^{\ast }\right) \in X\times X^{\ast }.$

\item The paire $\left( \left( x_{n},y_{n}^{\ast }\right) \right) _{n\in 
\mathbb{I}}$ of $X$ is said to be a besselian Schauder frame of $X$ if it is
both a Schauder frame and a besselian paire.
\end{enumerate}

\textbf{Remark. }For a besselian paire $\mathcal{F}:=\left( \left(
x_{n},y_{n}^{\ast }\right) \right) _{n\in \mathbb{I}}$ of $X$, the quantity 
\begin{equation*}
\mathcal{L}_{\mathcal{F}}:=\underset{\left( u,v^{\ast }\right) \in \mathbb{B}%
_{X}\times \mathbb{B}_{X^{\ast }}}{\sup }\left( \underset{n\in \mathbb{I}}{%
\sum }\left \vert y_{n}^{\ast }\left( u\right) \right \vert \left \vert
v^{\ast }\left( x_{n}\right) \right \vert \right)
\end{equation*}%
is finite and for each $\left( x,y\right) \in X\times X^{\ast }$, the
following inequality holds%
\begin{equation*}
\underset{n\in \mathbb{I}}{\sum }\left \vert y_{n}^{\ast }\left( x\right)
\right \vert \left \vert y^{\ast }\left( x_{n}\right) \right \vert \leq 
\mathcal{L}_{\mathcal{F}}\left \Vert x\right \Vert _{X}\left \Vert y^{\ast
}\right \Vert _{X^{\ast }}
\end{equation*}%
The quantity $\mathcal{L}_{\mathcal{F}}$ is then called the constant of the
besselian paire $\mathcal{F}$.

For all the material on Banach spaces or Hilbertian frames, one can refer to 
\cite{meg}, \cite{lin01}, \cite{lin02}, \cite{woj} and \cite{chr3}. In the
sequel $\left( E,\left \Vert \cdot \right \Vert _{E}\right) $ will represent
a given separable Banach space, $\left( \left( a_{n},b_{n}^{\ast }\right)
\right) _{n\in 
\mathbb{N}
^{\ast }}$a paire of $E$ and $p\in \left] 1,+\infty \right[ $ a given
constant. When $\left( \left( a_{n},b_{n}^{\ast }\right) \right) _{n\in 
\mathbb{N}
^{\ast }}$is a Schauder frame of $E,$ we will always denote, for each $n\in 
\mathbb{N}^{\ast },$ by $S_{n}$ the operator : 
\begin{equation*}
\begin{array}{cccc}
S_{n}: & E & \rightarrow & E \\ 
& x & \mapsto & \underset{j=1}{\overset{n}{\sum }}b_{j}^{\ast }\left(
x\right) a_{j}%
\end{array}%
\end{equation*}

\section{Examples of besselian Schauder frames}

\textbf{Example 3.1. }Let $E$ be a vector space over $\mathbb{K}$, and of
finite dimension $N$. Let $(f_{j})_{1\leq j\leq N}$ be a basis of $E$ and $%
(f_{j}^{\ast })_{1\leq j\leq N}$ its dual basis. Let $\Vert \cdot \Vert $ be
a norm on $E$. We denote by $\Vert \cdot \Vert _{1}$ (resp. $\Vert \cdot
\Vert _{1}^{\ast }$) the norm defined on $E$ (resp $E^{\ast }$), for each $%
\left( x_{1},...,x_{N}\right) \in \mathbb{K}^{N},$ by the formulas : 
\begin{equation*}
\left \Vert \underset{j=1}{\overset{N}{\sum }}x_{j}f_{j}\right \Vert _{1}=%
\underset{j=1}{\overset{N}{\sum }}\left \vert x_{j}\right \vert ,\left \Vert 
\underset{j=1}{\overset{N}{\sum }}x_{j}f_{j}^{\ast }\right \Vert _{1}^{\ast
}=\underset{j=1}{\overset{N}{\sum }}\left \vert x_{j}\right \vert _{N}
\end{equation*}%
For each $x=\underset{j=1}{\overset{N}{\sum }}x_{j}f_{j}$ and $y^{\ast }=%
\underset{j=1}{\overset{N}{\sum }}y_{j}f_{j}^{\ast }$, we have 
\begin{equation*}
x=\underset{j=1}{\overset{N}{\sum }}f_{j}^{\ast }(x)f_{j}
\end{equation*}%
and 
\begin{align*}
& \underset{j=1}{\overset{N}{\sum }}\left \vert f_{j}^{\ast }(x)\right \vert
\left \vert y^{\ast }(f_{j})\right \vert =\underset{j=1}{\overset{N}{\sum }}%
\left \vert x_{j}\right \vert \left \vert y_{j}\right \vert \\
& \leq \left( \underset{j=1}{\overset{N}{\sum }}\left \vert x_{j}\right
\vert \right) \left( \underset{j=1}{\overset{N}{\sum }}\left \vert
y_{j}\right \vert \right) \\
& \leq \Vert x\Vert _{1}\Vert y^{\ast }\Vert _{1}^{\ast }
\end{align*}%
Since $E$ and $E^{\ast }$ are finite dimensional, there exists a constant $%
M>0$ independent of $x$ and $y$ such that: $\Vert x\Vert _{1}\leq M\Vert
x\Vert _{E}$ and $\Vert y^{\ast }\Vert _{1}^{\ast }\leq M\Vert y^{\ast
}\Vert _{E^{\ast }}$. It follows that: 
\begin{equation*}
\underset{j=1}{\overset{N}{\sum }}\left \vert f_{j}^{\ast }(x)\right \vert
\left \vert y^{\ast }(f_{j})\right \vert \leq M^{2}\Vert x\Vert _{E}\Vert
y^{\ast }\Vert _{E^{\ast }}
\end{equation*}%
Consequently, $\left( \left( f_{j},f_{j}^{\ast }\right) \right) _{1\leq
j\leq N}$ is a besselian Schauder frame of $E$.

\textbf{Example 3.2. }Let us consider the well-known Banach space $%
l^{1}\left( \mathbb{K}\right) .$ It is well-known that the dual of $%
l^{1}\left( \mathbb{K}\right) $ is isometrically isomorphic to th Banach
space $l^{\infty }\left( \mathbb{K}\right) $. More precisely, the following
mapping :%
\begin{equation*}
\begin{array}{cccc}
\varphi : & l^{\infty }\left( \mathbb{K}\right) & \rightarrow & \left(
l^{1}\left( \mathbb{K}\right) \right) ^{\ast } \\ 
& \left( \alpha _{n}\right) _{n\in \mathbb{N}^{\ast }} & \mapsto & \varphi
\left( \left( \alpha _{n}\right) _{n\in \mathbb{N}^{\ast }}\right)%
\end{array}%
\end{equation*}%
defined for each $\left( x_{n}\right) _{n\in \mathbb{N}^{\ast }}\in
l^{1}\left( \mathbb{K}\right) $ by the formula 
\begin{equation*}
\varphi \left( \left( \alpha _{n}\right) _{n\in \mathbb{N}^{\ast }}\right)
\left( \left( x_{n}\right) _{n\in \mathbb{N}^{\ast }}\right) =\underset{n=1}{%
\overset{+\infty }{\sum }}\alpha _{n}x_{n}
\end{equation*}%
is an isometric isomorphism from $\left( l^{\infty }\left( \mathbb{K}\right)
,\left \Vert \cdot \right \Vert _{l^{\infty }\left( \mathbb{K}\right)
}\right) $ onto $\left( \left( l^{1}\left( \mathbb{K}\right) \right) ^{\ast
},\left \Vert \cdot \right \Vert _{\left( l^{1}\left( \mathbb{K}\right)
\right) ^{\ast }}\right) $. Let 
\begin{equation*}
\left( e_{n}\right) _{n\in \mathbb{N}^{\ast }}:=\left( \left( \delta
_{nm}\right) _{m\in \mathbb{N}^{\ast }}\right) _{n\in \mathbb{N}^{\ast }}
\end{equation*}%
be the canonical Schauder basis of $l^{1}\left( \mathbb{K}\right) .$ It is
clear that $\left( e_{n}\right) _{n\in \mathbb{N}^{\ast }}$ is a also
sequence of vectors of $l^{\infty }\left( \mathbb{K}\right) .$ We set, for
each $\ n\in \mathbb{N}^{\ast }$: 
\begin{equation*}
e_{n}^{\ast }:=\varphi \left( e_{n}\right)
\end{equation*}%
For each $y^{\ast }\in \left( l^{1}\left( \mathbb{K}\right) \right) ^{\ast }$
and $x=\left( x_{n}\right) _{n\in \mathbb{N}^{\ast }}\in l^{1}\left( \mathbb{%
K}\right) ,$ we have : 
\begin{equation*}
x=\underset{n=1}{\overset{+\infty }{\sum }}e_{n}^{\ast }\left( x\right) e_{n}
\end{equation*}%
Furthermore we have : 
\begin{eqnarray*}
&&\underset{n=1}{\overset{+\infty }{\sum }}\left \vert e_{n}^{\ast }\left(
x\right) \right \vert \left \vert y^{\ast }\left( e_{n}\right) \right \vert
\leq \underset{n=1}{\overset{+\infty }{\sum }}\left \vert x_{n}\right \vert
\left \Vert y^{\ast }\right \Vert _{\left( l^{1}\left( \mathbb{K}\right)
\right) ^{\ast }}\left \Vert e_{n}\right \Vert _{l^{1}\left( \mathbb{K}%
\right) } \\
&\leq &\left \Vert x\right \Vert _{l^{1}\left( \mathbb{K}\right) }\left
\Vert y^{\ast }\right \Vert _{\left( l^{1}\left( \mathbb{K}\right) \right)
^{\ast }}
\end{eqnarray*}%
It follows that $\left( \left( e_{n},e_{n}^{\ast }\right) \right) _{n\in 
\mathbb{N}}$ is a besselian Schauder frame of $l^{1}\left( 
\mathbb{C}
\right) $.

\textbf{Example 3.3. } We consider now the Banach space $l^{p}\left( \mathbb{%
K}\right) $. It is well-known that the dual of $l^{p}\left( \mathbb{K}%
\right) $ is isometrically isomorphic to th Banach space $l^{p^{\ast
}}\left( \mathbb{K}\right) $. More precisely, the following mapping :%
\begin{equation*}
\begin{array}{cccc}
\psi : & l^{p^{\ast }}\left( \mathbb{K}\right) & \rightarrow & \left(
l^{p}\left( \mathbb{K}\right) \right) ^{\ast } \\ 
& \left( \alpha _{n}\right) _{n\in \mathbb{N}^{\ast }} & \mapsto & \psi
\left( \left( \alpha _{n}\right) _{n\in \mathbb{N}^{\ast }}\right)%
\end{array}%
\end{equation*}%
defined for each $\left( x_{n}\right) _{n\in \mathbb{N}^{\ast }}\in
l^{p}\left( \mathbb{K}\right) $ by the formula : 
\begin{equation*}
\psi \left( \left( \alpha _{n}\right) _{n\in \mathbb{N}^{\ast }}\right)
\left( \left( x_{n}\right) _{n\in \mathbb{N}^{\ast }}\right) =\underset{n=1}{%
\overset{+\infty }{\sum }}\alpha _{n}x_{n}
\end{equation*}%
is a well-defined isometric isomorphism from $\left( l^{p^{\ast }}\left( 
\mathbb{K}\right) ,\left \Vert \cdot \right \Vert _{l^{p^{\ast }}\left( 
\mathbb{K}\right) }\right) $ onto $\left( \left( l^{p}\left( \mathbb{K}%
\right) \right) ^{\ast },\left \Vert \cdot \right \Vert _{\left( l^{p}\left( 
\mathbb{K}\right) \right) ^{\ast }}\right) $. The sequence : 
\begin{equation*}
\left( e_{n}\right) _{n\in \mathbb{N}^{\ast }}:=\left( \left( \delta
_{nm}\right) _{m\in \mathbb{N}^{\ast }}\right) _{n\in \mathbb{N}^{\ast }}
\end{equation*}%
is the canonical Schauder basis of $l^{p}\left( \mathbb{K}\right) .$ It is
clear that $\left( e_{n}\right) _{n\in \mathbb{N}^{\ast }}$ is a also the
canonical Schauder basis of $l^{p^{\ast }}\left( \mathbb{K}\right) .$ We
set, for each $n\in \mathbb{N}^{\ast }:$%
\begin{equation*}
u_{n}^{\ast }:=\psi \left( e_{n}\right)
\end{equation*}%
Then it is clear that we have for each $x\in l^{p}\left( \mathbb{K}\right) $%
, $y^{\ast }\in \left( l^{p}\left( \mathbb{K}\right) \right) ^{\ast }:$%
\begin{equation*}
\left \{ 
\begin{array}{c}
\left \Vert x\right \Vert _{l^{p}\left( \mathbb{K}\right) }=\left( \underset{%
n=1}{\overset{+\infty }{\sum }}\left \vert u_{n}^{\ast }\left( x\right)
\right \vert ^{p}\right) ^{\frac{1}{p}}\text{ } \\ 
\left \Vert y^{\ast }\right \Vert _{\left( l^{p}\left( \mathbb{K}\right)
\right) ^{\ast }}=\left( \underset{n=1}{\overset{+\infty }{\sum }}\left
\vert y^{\ast }\left( e_{n}\right) \right \vert ^{p^{\ast }}\right) ^{\frac{1%
}{p^{\ast }}}%
\end{array}%
\right.
\end{equation*}%
It follows then that : 
\begin{eqnarray*}
&&\underset{n=1}{\overset{+\infty }{\sum }}\left \vert u_{n}^{\ast }\left(
x\right) \right \vert \left \vert y^{\ast }\left( e_{n}\right) \right \vert
\leq \left( \underset{n=1}{\overset{+\infty }{\sum }}\left \vert u_{n}^{\ast
}\left( x\right) \right \vert ^{p}\right) ^{\frac{1}{p}}\left( \underset{n=1}%
{\overset{+\infty }{\sum }}\left \vert y^{\ast }\left( e_{n}\right) \right
\vert ^{p^{\ast }}\right) ^{\frac{1}{p^{\ast }}} \\
&\leq &\left \Vert x\right \Vert _{l^{p}\left( \mathbb{K}\right) }\left
\Vert y^{\ast }\right \Vert _{l^{p^{\ast }}\left( \mathbb{K}\right) }
\end{eqnarray*}%
\  \ It follows that $\left( \left( e_{n},u_{n}^{\ast }\right) \right) _{n\in 
\mathbb{N}^{\ast }}$ is a besselian Schauder frame of $l^{p}\left( \mathbb{K}%
\right) .$

\textbf{Example 3.4. }We consider now the Banach space $c_{0}\left( \mathbb{K%
}\right) $. It is well-known that the dual of $c_{0}\left( \mathbb{K}\right) 
$ is isometrically isomorphic to the Banach space $l^{1}\left( \mathbb{K}%
\right) $. More precisely, the following mapping :%
\begin{equation*}
\begin{array}{cccc}
\chi : & l^{1}\left( \mathbb{K}\right) & \rightarrow & \left( c_{0}\left( 
\mathbb{K}\right) \right) ^{\ast } \\ 
& \left( \alpha _{n}\right) _{n\in \mathbb{N}^{\ast }} & \mapsto & \chi
\left( \left( \alpha _{n}\right) _{n\in \mathbb{N}^{\ast }}\right)%
\end{array}%
\end{equation*}%
defined for each $\left( x_{n}\right) _{n\in \mathbb{N}^{\ast }}\in
l^{p}\left( \mathbb{K}\right) $ by the formula : 
\begin{equation*}
\chi \left( \left( \alpha _{n}\right) _{n\in \mathbb{N}^{\ast }}\right)
\left( \left( x_{n}\right) _{n\in \mathbb{N}^{\ast }}\right) =\underset{n=1}{%
\overset{+\infty }{\sum }}\alpha _{n}x_{n}
\end{equation*}%
is a well-defined isometric isomorphism from $\left( l^{1}\left( \mathbb{K}%
\right) ,\left \Vert \cdot \right \Vert _{l^{1}\left( \mathbb{K}\right)
}\right) $ onto $\left( \left( c_{0}\left( \mathbb{K}\right) \right) ^{\ast
},\left \Vert \cdot \right \Vert _{\left( c_{0}\left( \mathbb{K}\right)
\right) ^{\ast }}\right) $. The sequence : 
\begin{equation*}
\left( e_{n}\right) _{n\in \mathbb{N}^{\ast }}:=\left( \left( \delta
_{nm}\right) _{m\in \mathbb{N}^{\ast }}\right) _{n\in \mathbb{N}^{\ast }}
\end{equation*}%
is the canonical Schauder basis of $c_{0}\left( \mathbb{K}\right) $ and also
the canonical Schauder basis of $l^{1}\left( \mathbb{K}\right) .$ We set,
for each $n\in \mathbb{N}^{\ast }:$%
\begin{equation*}
v_{n}^{\ast }:=\chi \left( e_{n}\right)
\end{equation*}%
Then it is clear that we have for each $x\in c_{0}\left( \mathbb{K}\right) $%
, $y^{\ast }\in \left( c_{0}\left( \mathbb{K}\right) \right) ^{\ast }:$%
\begin{equation*}
\left \{ 
\begin{array}{c}
\left \Vert x\right \Vert _{c_{0}\left( \mathbb{K}\right) }=\underset{n\in 
\mathbb{N}
^{\ast }}{\sup }\left \vert v_{n}^{\ast }\left( x\right) \right \vert \text{ 
} \\ 
\left \Vert y^{\ast }\right \Vert _{\left( c_{0}\left( \mathbb{K}\right)
\right) ^{\ast }}=\underset{n=1}{\overset{+\infty }{\sum }}\left \vert
y^{\ast }\left( e_{n}\right) \right \vert \text{ }%
\end{array}%
\right.
\end{equation*}%
It follows then that : 
\begin{equation*}
\underset{n=1}{\overset{+\infty }{\sum }}\left \vert v_{n}^{\ast }\left(
x\right) \right \vert \left \vert y^{\ast }\left( e_{n}\right) \right \vert
\leq \left \Vert x\right \Vert _{c_{0}\left( \mathbb{K}\right) }\left \Vert
y^{\ast }\right \Vert _{\left( c_{0}\left( \mathbb{K}\right) \right) ^{\ast
}}
\end{equation*}%
It follows that\ $\left( \left( e_{n},v_{n}^{\ast }\right) \right) _{n\in 
\mathbb{N}}$ is a besselian Schauder frame of $c_{0}\left( \mathbb{K}\right)
.$

\textbf{Example 3.5. }Let $\left( \mathcal{H},\left \langle \cdot ,\cdot
\right \rangle \right) $ be a separable real or complex Hilbert space and $%
\mathcal{G}:=\left( x_{n}\right) _{n\in \mathbb{N}^{\ast }}$ a Hilbert frame
of $\mathcal{H}.$ That is $\left( x_{n}\right) _{n\in \mathbb{N}^{\ast }}$
is a sequence of vectors of $\mathcal{H}$ for which there exist two
constants $\alpha ,\beta >0$ such that the following condition holds for
every $x\in \mathcal{H}:$ 
\begin{equation*}
\alpha \left \Vert x\right \Vert _{\mathcal{H}}^{2}\leq \underset{n=1}{%
\overset{+\infty }{\sum }}\left \vert \left \langle x_{n},x\right \rangle
\right \vert ^{2}\leq \beta \left \Vert x\right \Vert _{\mathcal{H}}^{2}
\end{equation*}%
Then the operator 
\begin{equation*}
\begin{array}{cccc}
S_{\mathcal{G}}: & \mathcal{H} & \rightarrow & \mathcal{H} \\ 
& x & \mapsto & \underset{n=1}{\overset{+\infty }{\sum }}\left \langle
x_{n},x\right \rangle x_{n}%
\end{array}%
\end{equation*}%
is a well defined isomorphism of Banach spaces invertible which is an
autoadjoint operator of $\mathcal{H}$ (\cite{chr3}, pages 122-124). It
follows that we have for each $x\in \mathcal{H}:$ 
\begin{equation*}
x=\underset{n=1}{\overset{+\infty }{\sum }}\left \langle S_{\mathcal{G}%
}^{-1}\left( x_{n}\right) ,x\right \rangle x_{n}
\end{equation*}%
On the other hand, by virtue of the Cauchy-Schwarz inequality, we have for
each $\left( x,y\right) \in \mathcal{H\times H}$: 
\begin{align*}
& \underset{n=1}{\overset{+\infty }{\sum }}\left \vert T_{S_{\mathcal{G}%
}^{-1}\left( a_{n}\right) }\left( x\right) \right \vert \left \vert
T_{y}\left( a_{n}\right) \right \vert =\underset{n=1}{\overset{+\infty }{%
\sum }}\left \vert \left \langle S_{\mathcal{G}}^{-1}\left( a_{n}\right)
,x\right \rangle \right \vert \left \vert \left \langle a_{n},y\right
\rangle \right \vert \\
& \leq \left( \underset{n=1}{\overset{+\infty }{\sum }}\left \vert \left
\langle S_{\mathcal{G}}^{-1}\left( a_{n}\right) ,x\right \rangle \right
\vert ^{2}\right) ^{\frac{1}{2}}\left( \underset{n=1}{\overset{+\infty }{%
\sum }}\left \vert \left \langle a_{n},y\right \rangle \right \vert
^{2}\right) ^{\frac{1}{2}} \\
& \leq \left( \underset{n=1}{\overset{+\infty }{\sum }}\left \vert \left
\langle a_{n},S_{\mathcal{G}}^{-1}\left( x\right) \right \rangle \right
\vert ^{2}\right) ^{\frac{1}{2}}\left( \underset{n=1}{\overset{+\infty }{%
\sum }}\left \vert \left \langle a_{n},y\right \rangle \right \vert
^{2}\right) ^{\frac{1}{2}} \\
& \leq \alpha \left \Vert S_{\mathcal{G}}^{-1}\left( x\right) \right \Vert _{%
\mathcal{H}}\left \Vert y\right \Vert _{\mathcal{H}} \\
& \leq \alpha \left \Vert S_{\mathcal{G}}^{-1}\right \Vert _{\mathcal{H}%
}\left \Vert x\right \Vert _{\mathcal{H}}\left \Vert T_{y}\right \Vert _{%
\mathcal{H}^{\ast }}
\end{align*}%
Consequently $\left( \left( x_{n},T_{S_{\mathcal{G}}^{-1}\left( x_{n}\right)
}\right) \right) _{n\in \mathbb{N}^{\ast }}$ is a besselian Schauder frame
of $\mathcal{H}$.

\section{Fundamental results}

\textbf{Proposition 4.1.}\textit{\ Let }$\left( f_{n}\right) _{n\in \mathbb{N%
}^{\ast }}$\textit{\ and }$\left( g_{n}\right) _{n\in \mathbb{N}^{\ast }}$%
\textit{\ be a Bessel sequences of a Hilbert space} $\left( \mathcal{H}%
,\left \langle .,.\right \rangle \right) $. \textit{Then the paire }$\left(
\left( f_{n},T_{g_{n}}\right) \right) _{n\in \mathbb{N}^{\ast }}$\textit{\
is a besselian paire of} $\mathcal{H}$.

\begin{proof}
We have for each $x\in \mathcal{H}$ and $y\in \mathcal{H}$ : 
\begin{align*}
& \underset{n=1}{\overset{+\infty }{\sum }}\left \vert T_{g_{n}}\left(
x\right) \right \vert \left \vert T_{y}\left( f_{n}\right) \right \vert =%
\underset{n=1}{\overset{+\infty }{\sum }}\left \vert \left \langle
x,g_{n}\right \rangle \left \langle f_{n},y\right \rangle \right \vert \\
& \leq \sqrt{\underset{n=1}{\overset{+\infty }{\sum }}\left \vert \left
\langle x,g_{n}\right \rangle \right \vert ^{2}}\sqrt{\underset{n=1}{\overset%
{+\infty }{\sum }}\left \vert \left \langle f_{n},y\right \rangle \right
\vert ^{2}}
\end{align*}%
But there exists a constant $B>0$ such that: 
\begin{equation*}
\left \{ 
\begin{array}{l}
\underset{n=1}{\overset{+\infty }{\sum }}\left \vert \left \langle
x,g_{n}\right \rangle \right \vert ^{2}\leq B\left \Vert x\right \Vert _{%
\mathcal{H}}^{2},\; \;x\in \mathcal{H} \\ 
\underset{n=1}{\overset{+\infty }{\sum }}\left \vert \left \langle
f_{n},y\right \rangle \right \vert ^{2}\leq B\left \Vert y\right \Vert _{%
\mathcal{H}}^{2},\; \;y\in \mathcal{H}%
\end{array}%
\right.
\end{equation*}%
Consequently, we have for each $x\in \mathcal{H}$ and $y\in \mathcal{H}$: 
\begin{equation*}
\underset{n=1}{\overset{+\infty }{\sum }}\left \vert T_{g_{n}}\left(
x\right) \right \vert \left \vert T_{y}\left( f_{n}\right) \right \vert \leq
B\left \Vert x\right \Vert _{\mathcal{H}}\left \Vert y\right \Vert _{%
\mathcal{H}}
\end{equation*}%
Then $\left( \left( f_{n},T_{g_{n}}\right) \right) _{n\in \mathbb{N}^{\ast
}} $ is a besselian paire of $\mathcal{H}$.\newline
\end{proof}

We can generalize the last proposition to the Banach space $E$ by means of
the notions of weakly $p$-summable sequences of $E$ and *-weakly $p$%
-summable sequences of $E^{\ast }.$ We need first to prove the following
useful lemma.

\textbf{Lemma 4.2.}

\textit{Let }$\left( x_{n}\right) _{n\in 
\mathbb{N}
^{\ast }}$\textit{be a sequence of elements of \ }$E$ \textit{and} $\left(
f_{n}^{\ast }\right) _{n\in 
\mathbb{N}
^{\ast }}$\textit{be a sequence of elements of \ }$E^{\ast }$.\textit{\ }

1\textit{.\textbf{\ }The sequence }$\left( x_{n}\right) _{n\in 
\mathbb{N}
^{\ast }}$\textit{\ is weakly }$p$\textit{-summable in }$E$\textit{\ if and
only if there exists a constant }$C>0$ \textit{such that} :%
\begin{equation*}
\left( \overset{+\infty }{\underset{n=1}{\sum }}\left \vert f^{\ast }\left(
x_{n}\right) \right \vert ^{p}\right) ^{\frac{1}{p}}\leq C\left \Vert
f^{\ast }\right \Vert _{E^{\ast }},\text{ }f^{\ast }\in E^{\ast }
\end{equation*}

2.\textbf{\ }\textit{The sequence }$\left( f_{n}^{\ast }\right) _{n\in 
\mathbb{N}
^{\ast }}$\textit{\ is *-weakly }$p$\textit{-summable in }$E^{\ast }$\textit{%
\ if and only if there exists a constant }$C>0$ \textit{such that} :%
\begin{equation*}
\left( \overset{+\infty }{\underset{n=1}{\sum }}\left \vert f_{n}^{\ast
}\left( x\right) \right \vert ^{p}\right) ^{\frac{1}{p}}\leq C\left \Vert
x\right \Vert _{E},\text{ }x\in E
\end{equation*}

\begin{proof}
1. The sufficiency, in the first equivalence is obvious. We prove the
necessity by means of the closed graph theorem applied to the operator%
\begin{equation*}
\begin{array}{cccc}
U: & E^{\ast } & \rightarrow & l^{p}\left( \mathbb{K}\right) \\ 
& f^{\ast } & \mapsto & \left( f^{\ast }(x_{n})\right) _{n\in \mathbb{N}%
^{\ast }}%
\end{array}%
\end{equation*}
Indeed $U$ is well-defined on $E^{\ast }$ and is linear. Then let $\left(
f_{n}^{\ast }\right) _{n\in \mathbb{N}^{\ast }}$ be a sequence of elements
of $E^{\ast }$ such that $\left( f_{n}^{\ast }\right) _{n\in \mathbb{N}%
^{\ast }}$ is convergent in $E^{\ast }$ to $y^{\ast }$ and $\left(
U(f_{n}^{\ast })\right) _{n\in \mathbb{N}^{\ast }}$ is convergent in $%
l^{p}\left( \mathbb{K}\right) $ to $\left( \lambda _{j}\right) _{j\in 
\mathbb{N}^{\ast }}$. It follows that: 
\begin{equation*}
\underset{n\rightarrow +\infty }{\lim }\left \vert f_{n}^{\ast
}(x_{j})-\lambda _{j}\right \vert =0,\;j\in \mathbb{N}^{\ast }
\end{equation*}%
So $\underset{n\rightarrow +\infty }{\lim }f_{n}^{\ast }(x_{j})=\lambda _{j}$
for each $j\in \mathbb{N}^{\ast }$, but we have also: 
\begin{equation*}
\underset{n\rightarrow +\infty }{\lim }\left \Vert f_{n}^{\ast }-y^{\ast
}\right \Vert _{E^{\ast }}=0
\end{equation*}%
It follows that 
\begin{equation*}
y^{\ast }(x_{j})=\lambda _{j},\text{ }j\in \mathbb{N}^{\ast }
\end{equation*}%
Consequently, $\left( \lambda _{n}\right) _{n\in \mathbb{N}^{\ast }}$ $\in
l^{p}\left( \mathbb{K}\right) $ and $\left( U(f_{n}^{\ast })\right) _{n\in 
\mathbb{N}^{\ast }}$ is convergent in $E^{\ast }$ to $\left( y^{\ast
}(x_{j})\right) _{j\in \mathbb{N}^{\ast }}=U(y^{\ast })$. Hence the graph of 
$U$ is closed. It follows that the operator $U$ is continuous. That is,
there exists a constant $C>0$ such that: 
\begin{equation*}
\left \Vert U(f^{\ast })\right \Vert _{l^{p}}\leq C\left \Vert f^{\ast
}\right \Vert _{E^{\ast }},\; \;f^{\ast }\in E^{\ast }
\end{equation*}%
Consequently, 
\begin{equation*}
\left( \underset{n=1}{\overset{+\infty }{\sum }}\left \vert f^{\ast
}(x_{n})\right \vert ^{p}\right) ^{1/p}\leq C\left \Vert f^{\ast }\right
\Vert _{E^{\ast }},\; \;f^{\ast }\in E^{\ast }
\end{equation*}%
Hence we achieve the proof of the first equivalence of the proposition.

2. The sufficiency in the second equivalence is also obvious. In a similar
way to the proof of the first equivalence we prove the necessity by means of
the closed graph theorem applied to the operator%
\begin{equation*}
\begin{array}{cccc}
V: & E & \rightarrow & l^{p}\left( \mathbb{K}\right) \\ 
& x & \mapsto & \left( f_{n}^{\ast }(x)\right) _{n\in \mathbb{N}^{\ast }}%
\end{array}%
\end{equation*}
\end{proof}

\textbf{Proposition 4.3. \ }\textit{Let }$\left( x_{n}\right) _{n\in \mathbb{%
N}^{\ast }}$\textit{\ be sequence of elements of }$E$ \textit{which is
weakly }$p$\textit{-summable in }$E$\textit{\ and }$\left( y_{n}^{\ast
}\right) _{n\in \mathbb{N}^{\ast }}$\textit{\ be sequence of elements of }$%
E^{\ast }$ \textit{which is *-weakly }$p^{\ast }$\textit{-summable in }$%
E^{\ast }$\textit{. Then }$\left( \left( x_{n},y_{n}^{\ast }\right) \right)
_{n\in \mathbb{N}^{\ast }}$\textit{\ is a besselian paire of }$E$\textit{.}

\begin{proof}
Let $x\in E$ and $y^{\ast }\in E^{\ast }$. We have then by virtue of H\"{o}%
lder's inequality: 
\begin{equation*}
\underset{n=1}{\overset{+\infty }{\sum }}\left \vert y_{n}^{\ast }(x)\right
\vert \left \vert y^{\ast }(x_{n})\right \vert \leq \left( \underset{n=1}{%
\overset{+\infty }{\sum }}\left \vert y_{n}^{\ast }(x)\right \vert ^{p^{\ast
}}\right) ^{1/p^{\ast }}\left( \underset{n=1}{\overset{+\infty }{\sum }}%
\left \vert y^{\ast }(x_{n})\right \vert ^{p}\right) ^{1/p}
\end{equation*}%
Sience $\left( x_{n}\right) _{n\in \mathbb{N}^{\ast }}$ is weakly $p$%
-summable and $\left( y_{n}^{\ast }\right) _{n\in \mathbb{N}^{\ast }}$ is $%
\ast $-weakly $p^{\ast }$-summable, there exists a constant $C>0$ such that
\ the following inequalities hold for each $x\in E,$ $y^{\ast }\in E^{\ast }$
: 
\begin{equation*}
\left \{ 
\begin{array}{l}
\left( \underset{n=1}{\overset{+\infty }{\sum }}\left \vert y^{\ast
}(x_{n})\right \vert ^{p}\right) ^{1/p}\leq C\left \Vert y^{\ast }\right
\Vert _{E^{\ast }} \\ 
\left( \underset{n=1}{\overset{+\infty }{\sum }}\left \vert y_{n}^{\ast
}(x)\right \vert ^{p^{\ast }}\right) ^{1/p^{\ast }}\leq C\left \Vert x\right
\Vert _{E}%
\end{array}%
\right.
\end{equation*}%
It follows that: 
\begin{equation*}
\underset{n=1}{\overset{+\infty }{\sum }}\left \vert y_{n}^{\ast }(x)\right
\vert \left \vert y^{\ast }(x_{n})\right \vert \leq C^{2}\left \Vert x\right
\Vert _{E}\left \Vert y^{\ast }\right \Vert _{E^{\ast }}
\end{equation*}%
for each $x\in E$ and $y^{\ast }\in E^{\ast }$.\newline
Consequently, $\left( \left( x_{n},y_{n}^{\ast }\right) \right) _{n\in 
\mathbb{N}^{\ast }}$ is a besselian paire of $E$.\newline
\end{proof}

\textbf{Proposition 4.3. }\textit{Assume that }$E$ \textit{has a Schauder
frame}.\textit{\ Then }$E$\textit{\ has the approximation property.}

\begin{proof}
For each $n\in \mathbb{N}^{\ast }$ the operator $S_{n}$ is a finite rank
operator on $E$ and we have for every $x\in E$ 
\begin{equation*}
\left \Vert S_{n}(x)-x\right \Vert _{E}=\left \Vert \underset{j=n+1}{\overset%
{+\infty }{\sum }}b_{j}^{\ast }\left( x\right) a_{j}\right \Vert _{E}
\end{equation*}%
So the following relation holds for each $x\in E:$%
\begin{equation*}
\underset{n\rightarrow +\infty }{\lim }\left \Vert S_{n}(x)-x\right \Vert
_{E}=0
\end{equation*}%
It follows from \cite[Proposition 4.3,page:73]{rya} that $E$ has the
approximation property.\newline
\end{proof}

\textbf{Proposition 4.5. }\textit{Let }$E$ \textit{be a Banach space} 
\textit{and} $\left( \left( a_{n},b_{n}^{\ast }\right) \right) _{n\in 
\mathbb{N}
^{\ast }}$\textit{be a Schauder frame of }$E$. \textit{If the following
condition holds }%
\begin{equation}
\underset{n=1}{\overset{+\infty }{\sum }}\left \Vert a_{n}\right \Vert
_{E}\left \Vert b_{n}\right \Vert _{E^{\ast }}<+\infty  \label{C}
\end{equation}
\textit{then }$E$\textit{\ will be finite dimensional.}

\begin{proof}
For each $n\in 
\mathbb{N}
^{\ast }$, $S_{n}$ is a finite rank linear operator. On the other hand we
have for every $n\in 
\mathbb{N}
^{\ast }:$ 
\begin{equation*}
\left \Vert S_{n}-Id_{E}\right \Vert _{L\left( E\right) }\leq \underset{j=n+1%
}{\overset{+\infty }{\sum }}\left \Vert a_{j}\right \Vert _{E}\left \Vert
b_{j}\right \Vert _{E^{\ast }}
\end{equation*}%
It follows, according to the assumption (\ref{C}) of the proposition, that
the identity mapping $Id_{E}$ belongs to the topological closure in $L\left(
E\right) $ of the set of finite rank operators. Hence $Id_{E}$ is a compact
operator. Consequently the Banach space $E$ is finite dimensional.
\end{proof}

\textbf{Remark. }Using the concept of nuclear operators between Banach
spaces, we can give another version of the proof of the proposition 4.5.

\begin{proof}
Under the assumption (\ref{C}) of the proposition 4.5., the tensor product $%
\underset{n=1}{\overset{+\infty }{\sum }}b_{n}^{\ast }\otimes _{\pi }a_{n}$
is well-defined, belongs to $E^{\ast }\widehat{\otimes }_{\pi }E$ and
fullfiles the relation :%
\begin{equation*}
Id_{E}=J\left( \underset{n=1}{\overset{+\infty }{\sum }}b_{n}^{\ast }\otimes
_{\pi }a_{n}\right)
\end{equation*}%
where $J$ is the operator $J:E^{\ast }\widehat{\otimes }_{\pi }E\rightarrow
L\left( E\right) $ that associates with the tensor $u:=\underset{n=1}{%
\overset{+\infty }{\sum }}y_{n}^{\ast }\otimes _{\pi }x_{n}\in E^{\ast }%
\widehat{\otimes }_{\pi }E$ the linear operator $F_{u}$ defined by 
\begin{equation*}
F_{u}\left( x\right) :=\underset{n=1}{\overset{+\infty }{\sum }}y_{n}^{\ast
}\left( x\right) x_{n},\text{ }x\in E
\end{equation*}%
It follows that $Id_{E}$ is a nuclear operator of $L\left( E\right) $. But,
we know that all the nuclear operator of of $L\left( E\right) $ are compact (%
\cite{rya}, Chapter 2, page 42). It follows that $Id_{E}$ is compact.
Consequently $E$ is finite dimensional.
\end{proof}

\textbf{Proposition 4.6.} \textit{We assume that }$E$\textit{\ is a weakly
sequentially complete Banach space and that} $\mathcal{F}:=\left( \left(
a_{n},b_{n}^{\ast }\right) \right) _{n\in 
\mathbb{N}
^{\ast }}$ \textit{is a besselian paire of }$E.$\textit{Then for each }$x\in
E$\textit{, the series }$\sum b_{n}^{\ast }\left( x\right) a_{n}$\textit{\
is unconditionally convergent in }$E.$

\begin{proof}
For each $x\in E$, $y\in E^{\ast }$ we have: 
\begin{align*}
\overset{+\infty }{\underset{n=1}{\sum }}\left \vert y\left( b_{n}^{\ast
}\left( x\right) a_{n}\right) \right \vert & =\overset{+\infty }{\underset{%
k=1}{\sum }}\left \vert b_{n}^{\ast }\left( x\right) y\left( a_{n}\right)
\right \vert \\
& \leq \mathcal{L}_{\mathfrak{F}}\left \Vert x\right \Vert _{E}\left \Vert
y\right \Vert _{E^{\ast }} \\
& <+\infty
\end{align*}%
Hence the series $\sum b_{n}^{\ast }\left( x\right) a_{n}$ is weakly
unconditionally convergent. Then, since $E$ is weakly sequentially complete,
the well-known Orlicz's theorem $(1929)$ (\cite[Proposition.4 ,page:59]{woj}%
, page 66) entails that the series $\sum b_{n}^{\ast }\left( x\right) a_{n}$
is unconditionally unconditionally convergent.
\end{proof}

\textbf{Proposition 4.7.}\textit{\ Let }$E_{1},...,E_{N}$\textit{\ be a
closed subspaces of }$E$\textit{\ such that }%
\begin{equation}
E=E_{1}\oplus ...\oplus E_{N}  \label{S}
\end{equation}

1\textit{. Assume that }$E$\textit{\ has a Schauder frame }$($\textit{resp.
a besselian Schauder frame}$)$\textit{\ then the Banach space }$E_{j}$%
\textit{\ }$($\textit{as a closed subspace of }$E)$\textit{\ has, for each }$%
j\in \left \{ 1,...,N\right \} $,\textit{\ a Schauder frame }$($\textit{%
resp. a besselian Schauder frame}$)$\textit{.}

2\textit{. Assume that the Banach space }$E_{j}$\textit{\ }$($\textit{as a
closed subspace of }$E)$\textit{\ has, for each }$j\in \left \{
1,...,N\right \} $,\textit{\ a Schauder frame }$($\textit{resp. a besselian
Schauder frame}$)$\textit{. Then the space }$E$\textit{\ has a Schauder
frame }$($\textit{resp. a besselian Schauder frame}$)$\textit{.}

\begin{proof}
1. a. It follows from (\ref{S}), that there exists a continuous projections $%
\varrho _{j}:E\rightarrow E,$ $j\in \left \{ 1,...,N\right \} $ such that $%
\varrho _{j}$ is a projection onto $E_{j}$ for each $j\in \left \{
1,...,N\right \} $ and $\underset{j=1}{\overset{N}{\sum }}\varrho
_{j}=Id_{E} $. Let $\mathcal{F}:=\left( \left( a_{n},b_{n}^{\ast }\right)
\right) _{n\in \mathbb{I}}$ be a Schauder frame of $E$. Then we have for
each $x\in E$\ : 
\begin{equation*}
x=\underset{n\in \mathbb{I}}{\sum }b_{n}^{\ast }(x)a_{n}
\end{equation*}

If $x\in E_{j}$ for some $j\in \left \{ 1,...,N\right \} $, then since $%
\varrho _{j}$ is continuous, we have :%
\begin{eqnarray*}
x &=&\varrho _{j}\left( x\right) \\
&=&\underset{n\in \mathbb{I}}{\sum }b_{n}^{\ast }(x)\varrho _{j}\left(
a_{n}\right)
\end{eqnarray*}

But we have for each $j\in \left \{ 1,...,N\right \} $ and $n\in 
\mathbb{N}
^{\ast }:$ 
\begin{equation*}
\varrho _{j}\left( a_{n}\right) \in E_{j},\text{ }b_{n\left \vert
E_{j}\right. }^{\ast }\in E_{j}^{\ast }
\end{equation*}%
It follows that the paire $\left( \left( \varrho _{j}\left( a_{n}\right)
,b_{n\left \vert E_{j}\right. }^{\ast }\right) \right) _{n\in 
\mathbb{N}
^{\ast }}$ is a Schauder frame of $E_{j}$, for each $j\in \left \{
1,...,N\right \} .$ \  \  \  \ 1. b. We assume now that $\left( \left(
a_{n},b_{n}^{\ast }\right) \right) _{n\in \mathbb{I}}$ is a besselian
Schauder frame of $E$. Then we have for each $x\in E$ and $y^{\ast }\in
E^{\ast }$: 
\begin{equation*}
\underset{n\in \mathbb{I}}{\sum }\left \vert b_{n}^{\ast }\left( x\right)
\right \vert \left \vert y^{\ast }\left( a_{n}\right) \right \vert \leq 
\mathcal{L}_{\mathcal{F}}\left \Vert x\right \Vert _{E}\left \Vert y^{\ast
}\right \Vert _{E^{\ast }}
\end{equation*}%
When $x\in E_{j}$ and $y^{\ast }\in E_{j}^{\ast }$ for some $j\in \left \{
1,...,N\right \} $, then there exists, thanks to Hahn-Banach theorem a
linear form $Y^{\ast }\in E^{\ast }$ such that :%
\begin{equation*}
\left \{ 
\begin{array}{c}
Y_{\left \vert E_{j}\right. }^{\ast }=y^{\ast } \\ 
\left \Vert Y^{\ast }\right \Vert _{E^{\ast }}=\left \Vert y^{\ast }\right
\Vert _{E_{j}^{\ast }}%
\end{array}%
\right.
\end{equation*}

It follows that we have for each $x\in E_{j}$ and $y^{\ast }\in E_{j}^{\ast
} $ : 
\begin{eqnarray*}
&&\underset{n\in \mathbb{I}}{\sum }\left \vert \left( b_{n\left \vert
E_{j}\right. }^{\ast }\right) \left( x\right) \right \vert \left \vert
y^{\ast }\left( \varrho _{j}\left( a_{n}\right) \right) \right \vert =%
\underset{n\in \mathbb{I}}{\sum }\left \vert \left( b_{n}^{\ast }\right)
\left( x\right) \right \vert \left \vert Y^{\ast }\left( \varrho _{j}\left(
a_{n}\right) \right) \right \vert \\
&=&\underset{n\in \mathbb{I}}{\sum }\left \vert \left( b_{n}^{\ast }\right)
\left( x\right) \right \vert \left \vert \left( Y^{\ast }\circ \varrho
_{j}\right) \left( a_{n}\right) \right \vert \\
&\leq &\mathcal{L}_{\mathcal{F}}\left \Vert x\right \Vert _{E}\left \Vert
Y^{\ast }\circ \varrho _{j}\right \Vert _{E^{\ast }} \\
&\leq &\mathcal{L}_{\mathcal{F}}\left \Vert \varrho _{j}\right \Vert
_{L\left( E\right) }\left \Vert x\right \Vert _{E}\left \Vert Y^{\ast
}\right \Vert _{E^{\ast }} \\
&\leq &\mathcal{L}_{\mathcal{F}}\left \Vert \varrho _{j}\right \Vert
_{L\left( E\right) }\left \Vert x\right \Vert _{E_{j}}\left \Vert y^{\ast
}\right \Vert _{E_{j}^{\ast }}
\end{eqnarray*}%
But the paire $\mathcal{F}_{j}:=\left( \left( \varrho _{j}\left(
a_{n}\right) ,b_{n\left \vert E_{j}\right. }^{\ast }\right) \right) _{n\in 
\mathbb{N}
^{\ast }}$ is already Schauder frame of $E_{j}$, for each $j\in \left \{
1,...,N\right \} .$ Consequently, $\mathcal{F}_{j}$ is a besselian Schauder
frame of $E_{j}$, for each $j\in \left \{ 1,...,N\right \} .$

2. a. Assume that, for each $j\in \left \{ 1,...,N\right \} $, the subspace $%
E_{j}$ has a Schauder frame $\mathcal{F}_{j}:=\left( \left(
a_{n,j},b_{n,j}^{\ast }\right) \right) _{n\in \mathbb{I}_{j}}$ where $%
\mathbb{I}_{j}$ is a set of the form $\mathbb{I}_{j}=\left \{
1,...,r_{j}\right \} $ $\left( r_{j}\in 
\mathbb{N}
^{\ast }\right) $ or $\mathbb{I}_{j}=$ $%
\mathbb{N}
^{\ast }.$\ Let us set then for each $j\in \left \{ 1,...,N\right \} $ 
\begin{equation*}
\mathcal{R}_{j}:=\left( \left( x_{n,j},y_{n,j}^{\ast }\right) \right) _{n\in 
\mathbb{%
\mathbb{N}
}^{\ast }}
\end{equation*}%
where $:$ 
\begin{equation*}
\left( x_{n,j},y_{n,j}^{\ast }\right) =\left \{ 
\begin{array}{c}
\left( a_{n,j},b_{n,j}^{\ast }\right) \text{ if }n\in \mathbb{I}_{j} \\ 
\left( 0_{E},0_{E^{\ast }}\right) \text{ if }n\in \mathbb{%
\mathbb{N}
}^{\ast }\backslash \mathbb{I}_{j}%
\end{array}%
\right.
\end{equation*}%
Then $\mathcal{R}_{j}$ is, for each $j\in \left \{ 1,...,N\right \} $, is a
Schauder frame of $E_{J}.$

On theother hand we have for each $x\in E$ and $j\in \left \{
1,...,N\right
\} :$%
\begin{equation*}
\varrho _{j}\left( x\right) =\underset{n=1}{\overset{+\infty }{\sum }}%
b_{n,j}^{\ast }\left( \varrho _{j}\left( x\right) \right) a_{n,j}
\end{equation*}%
It follows that we have for each $x\in E$ : 
\begin{eqnarray*}
x &=&\underset{j=1}{\overset{N}{\sum }}\varrho _{j}\left( x\right) \\
&=&\underset{j=1}{\overset{N}{\sum }}\underset{n=1}{\overset{+\infty }{\sum }%
}y_{n,j}^{\ast }\left( \varrho _{j}\left( x\right) \right) x_{n,j}
\end{eqnarray*}%
Let us now consider the paire $\left( \left( x_{n},y_{n}^{\ast }\right)
\right) _{n\in 
\mathbb{N}
^{\ast }}$ defined by :%
\begin{equation*}
\left \{ 
\begin{array}{c}
x_{kN+l}=x_{k+1,l} \\ 
y_{kN+l}^{\ast }=y_{k+1,l}^{\ast }\circ \varrho _{l}%
\end{array}%
\right.
\end{equation*}%
for $k\in 
\mathbb{N}
$ and $l\in \left \{ 1,...,N\right \} .$ It is clear that the series $\sum
y_{n}^{\ast }\left( x\right) x_{n}$ is convergent for each $x\in E$ to the
vector $x.$ Hence the paire $\left( \left( x_{n},y_{n}^{\ast }\right)
\right) _{n\in 
\mathbb{N}
^{\ast }}$ is a Schauder frame of $E.$

2. b. Assume now that, for each $j\in \left \{ 1,...,N\right \} ,$ the paire 
$\mathcal{F}_{j}$ is a besselian Schauder frame of the subspace $E_{j}$.
Then $\mathcal{R}_{j}$ is, for each $j\in \left \{ 1,...,N\right \} ,$ a
besselian Schauder frame of the subspace $E_{j}.$ Hence we have for each $%
x\in E$ and $j\in \left \{ 1,...,N\right \} :$%
\begin{eqnarray*}
&&\underset{n=1}{\overset{+\infty }{\sum }}\left \vert b_{n,j}^{\ast }\left(
\varrho _{j}\left( x\right) \right) \right \vert \left \vert y^{\ast }\left(
a_{n,j}\right) \right \vert \leq \mathcal{L}_{\mathcal{F}_{j}}\left \Vert
\varrho _{j}\left( x\right) \right \Vert _{E_{j}}\left \Vert y_{\left \vert
E_{j}\right. }^{\ast }\right \Vert _{E_{j}^{\ast }} \\
&\leq &\mathcal{L}_{\mathcal{F}_{j}}\left \Vert \varrho _{j}\right \Vert
_{L\left( E\right) }\left \Vert x\right \Vert _{E}\left \Vert y^{\ast
}\right \Vert _{E^{\ast }}
\end{eqnarray*}%
It follows that we have for each $x\in E:$%
\begin{eqnarray*}
&&\underset{n=1}{\overset{+\infty }{\sum }}\left \vert b_{n}^{\ast }\left(
x\right) \right \vert \left \vert y^{\ast }\left( a_{n}\right) \right \vert =%
\underset{j=1}{\overset{N}{\sum }}\underset{n=1}{\overset{+\infty }{\sum }}%
\left \vert b_{n,j}^{\ast }\left( \varrho _{j}\left( x\right) \right) \right
\vert \left \vert y^{\ast }\left( a_{n,j}\right) \right \vert \\
&\leq &\left( \underset{j=1}{\overset{N}{\sum }}\mathcal{L}_{\mathcal{F}%
_{j}}\left \Vert \varrho _{j}\right \Vert _{L\left( E\right) }\right) \left
\Vert x\right \Vert _{E}\left \Vert y^{\ast }\right \Vert _{E^{\ast }}
\end{eqnarray*}%
It follow that is a besselian Schauder frame of $E.$
\end{proof}

\textbf{Corollary 4.1. }\textit{Let }$E$\textit{\ be a Banach space with a
Schauder frame }$($\textit{resp. a besselian Schauder frame}$)$\textit{\ and
a subspace }$F$\textit{\ of }$E$\textit{\ which is complemented in }$E$%
\textit{. Then }$F$\textit{\ has a Schauder frame }$($\textit{resp. a
besselian Schauder frame}$)$\textit{.}

\begin{proof}
The corollary is a direct consequence of the previous proposition (with $N=2$%
).
\end{proof}

\textbf{Corollary 4.2. }\textit{Let }$E_{1},...,E_{N}$ \textit{be a Banach
spaces. Then the Banach space }$E_{1}\times ...\times E_{N}$\textit{\
(endowed with the norm }$\left \Vert \left( x_{1},...,x_{N}\right)
\right
\Vert _{E_{1}\times ...\times E_{N}}:=\underset{j=1}{\overset{N}{%
\sum }}\left \Vert x_{j}\right \Vert _{E_{j}})$\textit{\ has a Schauder
frame (resp. a besselian Schauder frame) if and only if }$E_{j}$\textit{\
has, for each }$j\in \left \{ 1,...,N\right \} $\textit{\ a Schauder frame
(resp. a besselian Schauder frame).}

\begin{proof}
We set for each $j\in \{1,...,N\}$. 
\begin{equation*}
\tilde{E}_{j}=\left \{ 
\begin{array}{c}
(x_{1},...,x_{N})\in E_{1}\times ...\times E_{N}:\; \\ 
x_{k}=0_{E}\; \text{if}\;k\neq j\; \text{and}\;x_{j}\in E_{j}%
\end{array}%
\right \}
\end{equation*}%
$\  \  \  \  \  \  \  \  \  \  \  \  \  \  \  \  \  \  \  \  \  \  \  \  \  \  \  \  \  \  \  \  \  \  \  \  \  \
\  \ $

1. For each $j\in \left \{ 1,...,N\right \} ,$ $\tilde{E}_{j}$ has a
Schauder frame (resp. a besselian Schauder frame) if and only if $E_{j}$ has
a Schauder frame (resp. a besselian Schauder frame);

2. Since $E_{1}\times ...\times E_{N}=\tilde{E}_{1}\oplus ...\oplus \tilde{E}%
_{N},$ it follows from the proposition above that $E_{1}\times ...\times
E_{N}$ has a Schauder frame (resp. a besselian Schauder frame) if and only
if $\tilde{E}_{j}$ has for every $j\in \{1,...,N\}$ a Schauder frame (resp.
a besselian Schauder frame) if and only if $\tilde{E}_{j}$ has for every $%
j\in \{1,...,N\}$ a Schauder frame (resp. a besselian Schauder frame).
\end{proof}

\section{Conditions for the existence of Schauder frames and besselian
Schauder frames.}

\textbf{Theorem 5.1.}\textit{\ }

1. \textit{If }$E$\textit{\ has a Schauder frame then there exists a Banach
space }$Z$\textit{\ such that }$Z$\textit{\ has a Schauder basis and }$E$%
\textit{\ is isometrically isomorphic to a complemented subspace of }$Z$%
\textit{.}

2. \textit{If }$E$\textit{\ is isomorphic to a complemented subspace of a
Banach space }$Z$\textit{\ which}\  \textit{has a Schauder basis then }$E$%
\textit{\ has a Schauder frame.}

\begin{proof}
1. Let $\mathcal{F}:=\left( (a_{n},b_{n}^{\ast })\right) _{n\in \mathbb{I}}$
be a Schauder frame of $E.$

If the set $\left \{ n\in \mathbb{I}\text{ : }a_{n}\neq 0_{E}\text{ and }%
b_{n}^{\ast }\neq 0_{E^{\ast }}\right \} $ is finite then $E$\ is finite
dimensional. In this case $E$ has an algebraic basis which is also a
Schauder basis. Hence, in this case, the assertion of the direct part of the
theorem is true for $Z=E.$

Now we can assume without loss of the generality that $\mathbb{I=%
\mathbb{N}
}^{\ast }$and 
\begin{equation*}
a_{n}\neq 0_{E},\text{ }b_{n}^{\ast }\neq 0_{E^{\ast }},\text{ }n\in 
\mathbb{N}
^{\ast }
\end{equation*}%
Let $Z_{E}$ be the set of the sequences $z:=\left( b_{n}^{\ast
}(x_{n})a_{n}\right) _{n\in \mathbb{N}^{\ast }}$ such that $x_{n}\in E$ for
each $n\in \mathbb{N}^{\ast }$ and the series $\sum b_{n}^{\ast
}(x_{n})a_{n} $ is convergent. It is clear that $Z_{E}$ is a vector space
over $\mathbb{K}, $ when endowed with the usual operations of addition of
sequences and multiplication of sequences by scalars. We can easily show
that the mapping $\left \Vert \cdot \right \Vert _{Z_{E}}$ $%
:Z_{E}\rightarrow 
\mathbb{R}
$ defined for each $z:=\left( b_{n}^{\ast }(x_{n})a_{n}\right) _{n\in 
\mathbb{N}^{\ast }}\in Z_{E}$ by : 
\begin{equation*}
\left \Vert z\right \Vert _{Z_{E}}:=\underset{n\in 
\mathbb{N}
^{\ast }}{\sup }\left \Vert \underset{j=1}{\overset{n}{\sum }}b_{j}^{\ast
}\left( x_{j}\right) a_{j}\right \Vert _{E}
\end{equation*}%
is well-defined and is a norm on the vector space $Z_{E}.$

\textbf{Claim 1.}\textit{\ }$\left( Z_{E},\left \Vert \cdot \right \Vert
_{Z_{E}}\right) $\textit{\ is a Banach space.}

\textbf{Proof of the claim :}

Let $\left( z_{m}\right) _{m\in \mathbb{N}^{\ast }}:=\left( \left(
b_{n}^{\ast }(x_{n,m})a_{n}\right) _{n\in \mathbb{N}^{\ast }}\right) _{m\in 
\mathbb{N}^{\ast }}$ be a Cauchy sequence in $\left( Z_{E},\left \Vert \cdot
\right \Vert _{Z_{E}}\right) $. Hence there exists, for each $\varepsilon >0$%
, an integer $N_{\varepsilon }\in \mathbb{N}^{\ast }$ such that : 
\begin{equation*}
\left \Vert z_{m}-z_{s}\right \Vert _{Z_{E}}\leq \frac{\varepsilon }{4},\;
\;s\geq m\geq N_{\epsilon }
\end{equation*}%
which means that the following inequality holds for each $j,n,s,m$ $\in 
\mathbb{N}
^{\ast }$ such that $s\geq m\geq N_{\varepsilon }:$ 
\begin{equation}
\left \Vert \underset{j=1}{\overset{n}{\sum }}b_{j}^{\ast }\left(
x_{j,m}-x_{j,s}\right) a_{j}\right \Vert _{Z_{E}}\leq \frac{\varepsilon }{4}
\label{a}
\end{equation}%
Hence the sequence $\left( b_{j}^{\ast }\left( x_{j,m}\right) \right) _{m\in 
\mathbb{N}^{\ast }}$ are, for each $j\in 
\mathbb{N}
^{\ast },$ a Cauchy sequences in $E$. Hence the sequences $\left(
b_{j}^{\ast }\left( x_{j,m}\right) \right) _{m\in \mathbb{N}^{\ast }}$ are
convergent in $E.$ We set then : 
\begin{equation*}
l_{j}:=\underset{m\rightarrow +\infty }{\lim }b_{j}^{\ast }\left(
x_{j,m}\right) ,\text{ }j\in 
\mathbb{N}
^{\ast }
\end{equation*}%
Since $b_{j}^{\ast }\neq 0_{E^{\ast }}$ for each $j\in \mathbb{N}^{\ast }$,
it follows that there exists for each $j\in \mathbb{N}^{\ast }$, $x_{j}\in E$
such that 
\begin{equation*}
b_{j}^{\ast }(x_{j})=l_{j},\; \;j\in \mathbb{N}^{\ast }
\end{equation*}%
Tending $s$ to infinity in (\ref{a}), we obtain for each $j,n,$ $m$ $\in 
\mathbb{N}
^{\ast }$ such that $m\geq N_{\varepsilon }:$ 
\begin{equation*}
\underset{n\in 
\mathbb{N}
^{\ast }}{\sup }\left \Vert \underset{j=1}{\overset{n}{\sum }}b_{j}^{\ast
}\left( x_{j,m}\right) a_{j}-\underset{j=1}{\overset{n}{\sum }}b_{j}^{\ast
}\left( x_{j}\right) a_{j}\right \Vert _{E}\leq \dfrac{\varepsilon }{4}
\end{equation*}%
It follows that we have for each integers $n_{2}\geq n_{1}:$ 
\begin{equation*}
\left \Vert \underset{j=n_{1}}{\overset{n_{2}}{\sum }}b_{j}^{\ast }\left(
x_{j,N_{\varepsilon }}\right) a_{j}-\underset{j=n_{1}}{\overset{n_{2}}{\sum }%
}b_{j}^{\ast }\left( x_{j}\right) a_{j}\right \Vert _{E}\leq \dfrac{%
\varepsilon }{2}
\end{equation*}%
But since the series $\sum b_{j}^{\ast }\left( x_{n,N_{\varepsilon }}\right)
a_{j}$ is convergent, it follows that there exists $M_{\varepsilon }\in 
\mathbb{N}^{\ast }$ such that the following inequality holds for each
integers$\;n_{1},n_{2}$ such that $n_{2}\geq n_{1}\geq M_{\varepsilon }:$ 
\begin{equation*}
\left \Vert \underset{j=n_{1}}{\overset{n_{2}}{\sum }}b_{j}^{\ast }\left(
x_{j,N_{\varepsilon }}\right) a_{j}\right \Vert _{E}\leq \dfrac{\varepsilon 
}{2}
\end{equation*}%
It follows that 
\begin{equation*}
\left \Vert \underset{j=n_{1}}{\overset{n_{2}}{\sum }}b_{j}^{\ast }\left(
x_{j}\right) a_{j}\right \Vert _{E}\leq \varepsilon ,\text{ }n_{2}\geq
n_{1}\geq M_{\varepsilon }
\end{equation*}%
Consequently, the series $\sum b_{n}^{\ast }\left( x_{n}\right) a_{n}$ is
convergent in $E$. Let us set $z_{\infty }:=\underset{n=1}{\overset{+\infty }%
{\sum }}b_{n}^{\ast }\left( x_{n}\right) a_{n}$ then $u\in Z_{E}$ and 
\begin{equation*}
\left \Vert z_{m}-z_{\infty }\right \Vert _{Z_{E}}\leq \dfrac{\varepsilon }{4%
},\; \;m\geq N_{\varepsilon }
\end{equation*}%
Hence $\left( Z_{E},\left \Vert \cdot \right \Vert _{Z_{E}}\right) $ is a
Banach space.\newline

$\square $

\textbf{Claim 2.} \textit{Let us set for each }$n\in N^{\ast }$\textit{: }%
\begin{equation*}
A_{n}:=\left( \delta _{nm}a_{n}\right) _{m\in \mathbb{N}^{\ast }}
\end{equation*}%
\textit{Then} $A_{n}\in Z_{E}$\textit{\ for each }$n\in 
\mathbb{N}
^{\ast }$\textit{\ and the sequence }$(A_{n})_{n\in \mathbb{N}^{\ast }}$%
\textit{\ is a monotone Schauder basis of} $Z_{E}$.

\textbf{Proof of the claim.}

Since $b_{n}^{\ast }\neq 0_{E^{\ast }}$ for each $n\in 
\mathbb{N}
^{\ast },$ there exists $y_{n}\in E$ such that $b_{n}^{\ast }(y_{n})=1$. We
set $m\in \mathbb{N}^{\ast },$ $c_{m}=\delta _{m,n}y_{n}$. It follows that 
\begin{equation*}
\left \{ 
\begin{array}{c}
A_{n}=\left( b_{m}^{\ast }(c_{m})a_{m}\right) _{m\in 
\mathbb{N}
^{\ast }} \\ 
\overset{m}{\underset{j=1}{\sum }}b_{j}^{\ast }(c_{j})a_{j}=a_{n},\text{ }%
m\geq n%
\end{array}%
\right.
\end{equation*}%
Consequently : 
\begin{equation*}
A_{n}\in Z_{E},\;n\in \mathbb{N}^{\ast }
\end{equation*}%
Let $z:=\left( b_{n}^{\ast }(x_{n})a_{n}\right) _{n\in \mathbb{N}^{\ast }}$ $%
\ $be an element of $Z_{E}$. Then we have for each $n\in \mathbb{N}^{\ast }$%
: 
\begin{equation*}
z-\underset{j=1}{\overset{n}{\sum }}b_{j}^{\ast }\left( x_{j}\right)
A_{j}=\left( 
\begin{array}{c}
\underset{n\; \text{times}}{\underbrace{0_{E},...,0_{E}}},b_{n+1}^{\ast
}(x_{n+1})a_{n+1}, \\ 
b_{n+2}^{\ast }(x_{n+2})a_{n+2},...%
\end{array}%
\right)
\end{equation*}%
It follows that: 
\begin{equation*}
\left \Vert {z-\underset{j=1}{\overset{n}{\sum }}b_{j}^{\ast }\left(
x_{j}\right) A_{j}}\right \Vert _{Z_{E}}=\underset{r\geq n+1}{\sup }\left
\Vert \underset{j=n+1}{\overset{r}{\sum }}b_{j}^{\ast }\left( x_{j}\right)
a_{j}\right \Vert _{E}
\end{equation*}%
It follows that we have in $Z_{E}$ 
\begin{equation*}
z=\underset{n=1}{\overset{+\infty }{\sum }}b_{n}^{\ast }\left( x_{n}\right)
A_{n}
\end{equation*}%
On the other hand, we have for each $(\alpha _{1},\alpha _{2},...,\alpha
_{n+1})\in \mathbb{K}^{n+1}$ $(n\in \mathbb{N}^{\ast })$ : 
\begin{align*}
& \left \Vert {\underset{s=1}{\overset{n}{\sum }}\alpha _{s}A_{s}}\right
\Vert _{Z_{E}}=\left \Vert (b_{1}^{\ast }(\alpha
_{1}y_{1})a_{1},...,b_{n}^{\ast }(\alpha
_{n}y_{n})a_{n},0_{E},0_{E},...)\right \Vert _{Z_{E}} \\
& =\underset{1\leq r\leq n}{\sup }\left \Vert \underset{j=1}{\overset{r}{%
\sum }}b_{j}^{\ast }\left( \alpha _{j}y_{j}\right) a_{j}\right \Vert _{E} \\
& \leq \underset{1\leq r\leq n+1}{\sup }\left \Vert \underset{j=1}{\overset{r%
}{\sum }}b_{j}^{\ast }\left( \alpha _{j}y_{j}\right) a_{j}\right \Vert _{E}
\\
& \leq \left \Vert \underset{j=1}{\overset{n+1}{\sum }}\alpha
_{j}a_{j}\right \Vert _{Z_{E}}
\end{align*}%
Consequently $(A_{n})_{n\in \mathbb{N}^{\ast }}$ is a monotone Schauder
basis of $Z_{E}$.

$\square $

\textbf{Claim} \textbf{3.}

1. \textit{The mapping }$:$%
\begin{equation*}
\begin{array}{cccc}
T_{0}: & E & \rightarrow & Z_{E} \\ 
& x & \mapsto & \left( b_{n}^{\ast }(x)a_{n}\right) _{n\in \mathbb{N}^{\ast
}}%
\end{array}%
\end{equation*}%
\textit{\ is well-defined with a closed image }$T_{0}\left( E\right) $%
\textit{. Furthermore }$T_{0}$\textit{\ is an isomorphism from }$E$\textit{\
onto }$T_{0}(E).$\textit{\ }

2. \textit{The subspace }$T_{0}(E)$\textit{\ is complemented in }$Z_{E}$%
\textit{.}

\textbf{Proof of the claim. }

1. We have for each $x\in E$ :%
\begin{equation}
\underset{n\rightarrow +\infty }{\lim }\left \Vert x-S_{n}\left( x\right)
\right \Vert _{E}=0  \label{d}
\end{equation}%
It follows, thanks to (\cite{pel3}, page $239$), that the quantity 
\begin{equation*}
\mathcal{K}_{\mathcal{F}}:=\underset{n\in 
\mathbb{N}
^{\ast }}{\sup }\left \Vert S_{n}\right \Vert _{L\left( E\right) }
\end{equation*}%
is finite, that is :%
\begin{equation*}
\mathcal{K}_{\mathcal{F}}:=\underset{n\in 
\mathbb{N}
^{\ast },\text{ }y\in \mathbb{B}_{E}}{\sup }\left \Vert \underset{j=1}{%
\overset{n}{\sum }}b_{j}^{\ast }\left( y\right) a_{j}\right \Vert
_{E}<+\infty
\end{equation*}%
It follows that we have for each $x\in E$ $:$ \textbf{\ }%
\begin{equation*}
\underset{n\in 
\mathbb{N}
^{\ast }}{\sup }\left \Vert \underset{j=1}{\overset{n}{\sum }}b_{j}^{\ast
}\left( x\right) a_{j}\right \Vert _{E}\leq \mathcal{K}_{\mathcal{F}}\left
\Vert x\right \Vert _{E}
\end{equation*}

Hence, for each $x\in E$, we have $\left( b_{n}^{\ast }(x)a_{n}\right)
_{n\in \mathbb{N}^{\ast }}\in Z_{E}.$ It follows that the mapping $T_{0}$ is
a well-defined linear operator of $E$ and fulfilles for each $x\in E$ the
inequality $:$%
\begin{eqnarray*}
\left \Vert T_{0}(x)\right \Vert _{Z_{E}} &=&\underset{n\in 
\mathbb{N}
^{\ast }}{\sup }\left \Vert \underset{j=1}{\overset{n}{\sum }}b_{j}^{\ast
}\left( x\right) a_{j}\right \Vert _{E} \\
&\leq &\mathcal{K}_{\mathcal{F}}\left \Vert x\right \Vert _{E}
\end{eqnarray*}%
On the other hand, we have for each $x\in E$ : 
\begin{equation*}
\left \Vert \underset{j=1}{\overset{n}{\sum }}b_{j}^{\ast }\left( x\right)
a_{j}\right \Vert _{E}\leq \left \Vert T_{0}(x)\right \Vert _{Z_{E}},\text{ }%
n\in 
\mathbb{N}
^{\ast }
\end{equation*}

It follows, in view of (\ref{d}), that :%
\begin{equation*}
\left \Vert x\right \Vert _{E}\leq \left \Vert T_{0}(x)\right \Vert _{Z_{E}},%
\text{ }x\in E
\end{equation*}%
Consequently we have :%
\begin{equation}
\left \Vert x\right \Vert _{E}\leq \left \Vert T_{0}(x)\right \Vert
_{Z_{E}}\leq \mathcal{K}_{\mathcal{F}}\left \Vert x\right \Vert _{E},\text{ }%
x\in E  \label{e}
\end{equation}%
\  \textit{\ }It follows, from (\ref{e}), that $T_{0}\left( E\right) $ is a
closed subspace of $E$ and that $T_{0}$ is an isomorphism from $E$ onto $%
T_{0}\left( E\right) $.

2. Let us consider the mapping $\varrho :Z_{E}\longrightarrow Z_{E}$ defined
for each $z:=\left( b_{n}^{\ast }(x_{n})a_{n}\right) _{n\in \mathbb{N}^{\ast
}}\in Z_{E}$ by the formula :%
\begin{equation*}
\varrho (z):=\left( b_{n}^{\ast }\left( \underset{m=1}{\overset{+\infty }{%
\sum }}b_{m}^{\ast }(x_{m})a_{m}\right) a_{n}\right) _{n\in \mathbb{N}^{\ast
}}
\end{equation*}%
It is easy to check that $\varrho $ is a well-defined linear mapping and
that : 
\begin{equation*}
\left \{ 
\begin{array}{c}
\varrho (z)=T_{0}\left( \left( \underset{m=1}{\overset{+\infty }{\sum }}%
b_{m}^{\ast }(x_{m})a_{m}\right) \right) \in T_{0}\left( E\right) ,\text{ }%
z\in Z_{E} \\ 
\varrho (z)=z\text{, }z\in T_{0}(E)%
\end{array}%
\right.
\end{equation*}%
It follows that $\varrho (Z_{E})=T_{0}(E).$ For each $z\in Z_{E}$, we have
also : 
\begin{eqnarray*}
&&\left \Vert \varrho (z)\right \Vert _{Z_{E}}=\left \Vert T_{0}\left(
\left( \underset{m=1}{\overset{+\infty }{\sum }}b_{m}^{\ast
}(x_{m})a_{m}\right) \right) \right \Vert _{Z_{E}} \\
&\leq &\mathcal{K}_{\mathcal{F}}\left \Vert \underset{m=1}{\overset{+\infty }%
{\sum }}b_{m}^{\ast }(x_{m})a_{m}\right \Vert _{Z_{E}} \\
&\leq &\mathcal{K}_{\mathcal{F}}\underset{n\in 
\mathbb{N}
^{\ast }}{\sup }\left \Vert \underset{j=1}{\overset{n}{\sum }}b_{j}^{\ast
}(x_{j})a_{j}\right \Vert _{E} \\
&\leq &\mathcal{K}_{\mathcal{F}}\left \Vert z\right \Vert _{Z_{E}}
\end{eqnarray*}%
It follows that $\varrho :Z_{E}\rightarrow Z_{E}$ is a bounded linear
operator.

On the other hand we have for each $z\in Z_{E}:$%
\begin{eqnarray*}
&&\left( \varrho \circ \varrho \right) \left( z\right) =\varrho \left(
\varrho \left( z\right) \right) \\
&=&\varrho \left( \left( b_{n}^{\ast }\left( \underset{m=1}{\overset{+\infty 
}{\sum }}b_{m}^{\ast }(x_{m})a_{m}\right) a_{n}\right) _{n\in \mathbb{N}%
^{\ast }}\right) \\
&=&\left( b_{n}^{\ast }\left( \underset{k=1}{\overset{+\infty }{\sum }}%
b_{k}^{\ast }\left( \underset{l=1}{\overset{+\infty }{\sum }}b_{l}^{\ast
}(x_{l})a_{l}\right) a_{k}\right) a_{n}\right) _{n\in \mathbb{N}^{\ast }} \\
&=&\left( b_{n}^{\ast }\left( \underset{l=1}{\overset{+\infty }{\sum }}%
b_{l}^{\ast }(x_{l})a_{l}\right) a_{n}\right) _{n\in \mathbb{N}^{\ast }} \\
&=&\varrho \left( z\right)
\end{eqnarray*}%
It follows that $\varrho \circ \varrho =\varrho .$

Consequenquently $\varrho $ is a bounded projection of $Z_{E}$ onto $%
T_{0}(E) $. It follows that $T_{0}(E)$ is a complemented subspace of $Z_{E}$.

$\square $

\textbf{End of the proof of the direct part of the theorem :}

In the claims above, it is proved that the Banach space $E$ is isomorphic to
the complemented subspace $T_{0}\left( E\right) $ of the Banach space $Z_{E}$
which has a monotone Schauder basis $\left( A_{n}\right) _{n\in 
\mathbb{N}
^{\ast }}$. But there exists, according to (\cite{pel1}, page 211,
proposition 1), a Banach space $Z$\ and an isomorphism $T_{1}:Z_{E}%
\rightarrow Z$ such that $E$ is isometrically isomorphic to the subspace $%
T_{1}\left( T_{0}\left( E\right) \right) .$ It follows that $\left(
T_{1}\left( A_{n}\right) \right) _{n\in 
\mathbb{N}
^{\ast }}$ is a Schauder basis of the Banach space $Z$ and that $T_{1}\left(
T_{0}\left( E\right) \right) $ is a complemented subspace of $Z.$ Hence we
achieve the proof of the direct part of the theorem.

2. We assume now that $E$\ is isomorphic to the complemented closed subspace 
$Y$ of a Banach space $Z$\ which\ has a Schauder basis $\mathcal{F}$:=$%
\left( \left( u_{n},v_{n}^{\ast }\right) \right) _{n\in \mathbb{I}}.$ Then $%
\mathcal{F}$\ is a Schauder frame of $Z$. It follows, thanks to the
corollary 4.1., that the Banach space $E$ has a Schauder frame.
\end{proof}

\textbf{Theorem 5.2. }\textit{Assume that }$E$\textit{\ is a weakly
sequentially complete Banach space such that }$E$ \textit{has a besselian
Schauder frame, then} $E$\textit{\ is isometrically isomorphic to a} \textit{%
complemented subspace of a Banach space }$W$\textit{\ wich has an
unconditional Schauder basis.}

\begin{proof}
Let $\mathcal{F}:=\left( (a_{n},b_{n}^{\ast })\right) _{n\in \mathbb{I}}$ be
a besselian Schauder frame of $E.$

If the set $\left \{ n\in \mathbb{I}\text{ : }a_{n}\neq 0_{E}\text{ and }%
b_{n}^{\ast }\neq 0_{E^{\ast }}\right \} $ \ is finite then $E$\ is finite
dimensional. In this case $E$ has an unconditional Schauder. Hence, in this
case, the assertion of the theorem is true for $W=E.$

Now we can assume without loss of the generality that $\mathbb{I=%
\mathbb{N}
}^{\ast }$and 
\begin{equation*}
a_{n}\neq 0_{E},\text{ }b_{n}^{\ast }\neq 0_{E^{\ast }},\text{ }n\in 
\mathbb{N}
^{\ast }
\end{equation*}%
We denote by $W_{E}$ the sequences $(b^{\ast }\left( x_{n}\right)
a_{n})_{n\in 
\mathbb{N}
^{\ast }}$ such that $x_{n}\in E$ for each $n\in 
\mathbb{N}
^{\ast }$ and the series $\sum b^{\ast }\left( x_{n}\right) a_{n}$ is
unconditionally convergent in $E.$ It is clear that $W_{E}$ is a vector
space over $\mathbb{K},$ when endowed with the usual operations of addition
of sequences and multiplication of sequences by scalars.

\textbf{Claim 1}: \textit{For each }$z=\left( b_{n}^{\ast
}(x_{n})a_{n}\right) _{n\in \mathbb{N}^{\ast }}$\textit{\ in }$W_{E}$\textit{%
, the quantity }%
\begin{equation*}
\Vert z\Vert _{W_{E}}:=\underset{\sigma \in \mathfrak{S}_{\mathbb{N}^{\ast
}},\;n\in \mathbb{N}^{\ast }}{\sup }\left \Vert \underset{j=1}{\overset{n}{%
\sum }}b_{\sigma (j)}^{\ast }\left( x_{\sigma (j)}\right) a_{\sigma
(j)}\right \Vert _{E}
\end{equation*}%
\textit{is finite.\newline
}\textbf{Proof of the claim :}

Since the series $\sum b_{n}^{\ast }(x_{n})a_{n}$ is unconditionally
convergent in $E$, it follows that there exists for each $\varepsilon >0$ an
integer $k(\varepsilon )\in \mathbb{N}^{\ast }$ such that the inegality 
\begin{equation*}
\left \Vert \sum_{j\in L}b_{j}^{\ast }(x_{j})a_{j}\right \Vert _{E}\leq
\varepsilon
\end{equation*}%
holds for each $L\in \mathcal{D}_{\mathbb{N}^{\ast }}$ which fulfilles the
condition $\min \left( L\right) \geq k(\varepsilon )+1$. Let $A\in \mathcal{D%
}_{\mathbb{N}^{\ast }}$ . Hence for $\varepsilon =1$, there exists $k(1)\in 
\mathbb{N}^{\ast }$ such that 
\begin{equation*}
\left \{ 
\begin{array}{l}
\left \Vert \underset{j\in A\cap \left \{ 1,...k(1)\right \} }{\sum }%
b_{j}^{\ast }(x_{j})a_{j}\right \Vert _{E}\leq \overset{k(1)}{\underset{j=1}{%
\sum }}\left \Vert b_{j}^{\ast }(x_{j})a_{j}\right \Vert _{E} \\ 
\left \Vert \underset{j\in A\backslash \left \{ 1,...k(1)\right \} }{\sum }%
b_{j}^{\ast }(x_{j})a_{j}\right \Vert _{E}\leq 1%
\end{array}%
\right.
\end{equation*}%
\newline
It follows that: 
\begin{equation*}
\left \Vert \underset{j\in A}{\sum }b_{j}^{\ast }(x_{j})a_{j}\right \Vert
_{E}\leq 1+\overset{k(1)}{\underset{j=1}{\sum }}\left \Vert b_{j}^{\ast
}(x_{j})a_{j}\right \Vert _{E}
\end{equation*}%
Consequently: 
\begin{equation*}
\underset{A\in \mathcal{D}_{\mathbb{N}^{\ast }}}{\sup }\left \Vert
\sum_{j\in A}b_{j}^{\ast }\left( x_{j}\right) a_{j}\right \Vert _{E}<+\infty
\end{equation*}%
It follows that: 
\begin{equation*}
\underset{\sigma \in \mathfrak{S}_{\mathbb{N}^{\ast }},\;n\in \mathbb{N}%
^{\ast }}{\sup }\left \Vert \underset{j=1}{\overset{n}{\sum }}b_{\sigma
(j)}^{\ast }\left( x_{\sigma (j)}\right) a_{\sigma (j)}\right \Vert
_{E}<+\infty
\end{equation*}%
Hence the quantity $\Vert z\Vert _{W_{E}}$ is finite.\newline

$\square $

\textbf{Claim 2} : \textit{The mapping }%
\begin{equation*}
\begin{array}{cccc}
\Vert .\Vert _{W_{E}}: & W_{E} & \rightarrow & \mathbb{R} \\ 
& z & \mapsto & \Vert z\Vert _{W_{E}}%
\end{array}%
\end{equation*}%
\textit{is a norm on }$W_{E}$ \textit{and }$(W_{E},\Vert .\Vert _{W_{E}})$%
\textit{\ is a Banach space.}

\textbf{Proof of the \ claim :}

Direct computations show that $\Vert .\Vert _{W_{E}}$ is a norm on $W_{E}$.
Let us show that $(W_{E},\Vert .\Vert _{W_{E}})$ is a Banach space. Indeed
let $\left( z_{m}\right) _{m\in 
\mathbb{N}
^{\ast }}:=\left( \left( b_{n}^{\ast }(x_{n,m})a_{n}\right) _{n\in \mathbb{N}%
^{\ast }}\right) _{m\in \mathbb{N}^{\ast }}$ be a Cauchy sequence in $%
(W_{E},\Vert .\Vert _{W_{E}})$. Then for each $\varepsilon >0$, there exists 
$N_{\varepsilon }\in \mathbb{N}^{\ast }$ such that the following inequality
holds for each $\sigma \in \mathfrak{S}_{\mathbb{N}^{\ast }}$, $s\geq r\geq
N_{\varepsilon }$ and $n\in \mathbb{N}^{\ast }$: 
\begin{equation*}
\left \Vert \underset{j=1}{\overset{n}{\sum }}\left( b_{\sigma (j)}^{\ast
}\left( x_{\sigma (j),s}\right) a_{\sigma (j)}-b_{\sigma (j)}^{\ast }\left(
x_{\sigma (j),r}\right) a_{\sigma (j)}\right) \right \Vert _{E}\leq \dfrac{%
\varepsilon }{4}
\end{equation*}%
It follows that: 
\begin{equation*}
\left \Vert b_{\sigma (n)}^{\ast }\left( x_{\sigma (n),s}\right) a_{\sigma
(n)}-b_{\sigma (n)}^{\ast }\left( x_{\sigma (n),r}\right) a_{\sigma
(n)}\right \Vert _{E}\leq \dfrac{\varepsilon }{2}
\end{equation*}%
for $s\geq r\geq N_{\varepsilon }$,$n\in \mathbb{N}^{\ast }$ and $\sigma \in 
\mathfrak{S}_{\mathbb{N}^{\ast }}$. Consequently, for each $n\in \mathbb{N}%
^{\ast }$ the sequence $\left( b_{n}^{\ast }\left( x_{n,r}\right)
a_{n}\right) _{n\in \mathbb{N}^{\ast }}$ is convergent in $E$ to a vector of
the form $b_{n}^{\ast }\left( x_{n}\right) a_{n}$. Tending $s$ to infinity,
we obtain for each $n,r\in \mathbb{N}^{\ast }$, $\sigma \in \mathfrak{S}_{%
\mathbb{N}^{\ast }}$ $\ $with $r\geq N_{\varepsilon },$ the following
inequality : 
\begin{equation*}
\left \Vert \underset{j=1}{\overset{n}{\sum }}\left( b_{\sigma (j)}^{\ast
}\left( x_{\sigma (j)}\right) a_{\sigma (j)}-b_{\sigma (j)}^{\ast }\left(
x_{\sigma (j),r}\right) a_{\sigma (j)}\right) \right \Vert _{E}\leq \dfrac{%
\varepsilon }{4}
\end{equation*}%
It follows, for each $n_{1},n_{2}\in 
\mathbb{N}
^{\ast },$ $n_{2}\geq n_{1}$, that : 
\begin{equation*}
\left \Vert \underset{j=n_{1}}{\overset{n_{2}}{\sum }}b_{\sigma (j)}^{\ast
}\left( x_{\sigma (j)}\right) a_{\sigma (j)}-\underset{j=n_{1}}{\overset{%
n_{2}}{\sum }}b_{\sigma (j)}^{\ast }\left( x_{\sigma (j),N_{\varepsilon
}}\right) a_{\sigma (j)}\right \Vert _{E}\leq \dfrac{\varepsilon }{2}
\end{equation*}%
But the series $\sum_{j}b_{\sigma (j)}^{\ast }\left( x_{\sigma
(j),N_{\varepsilon }}\right) a_{\sigma (j)}$ is convergent for each $\sigma
\in \mathfrak{S}_{\mathbb{N}^{\ast }}$, hence there exists for each $\sigma
\in \mathfrak{S}_{\mathbb{N}^{\ast }},$ an integer $M_{\varepsilon ,\sigma
}\in \mathbb{N}^{\ast }$ such that: 
\begin{equation*}
\left \Vert \underset{j=n_{1}}{\overset{n_{2}}{\sum }}b_{\sigma (j)}^{\ast
}\left( x_{\sigma (j),N_{\varepsilon }}\right) a_{\sigma (j)}\right \Vert
_{E}\leq \dfrac{\varepsilon }{2},\;n_{2}\geq n_{1}\geq M_{\varepsilon
,\sigma }
\end{equation*}%
It follows that: 
\begin{equation*}
\left \Vert \underset{j=n_{1}}{\overset{n_{2}}{\sum }}b_{\sigma (j)}^{\ast
}\left( x_{\sigma (j)}\right) a_{\sigma (j)}\right \Vert _{E}\leq
\varepsilon ,\text{ }n_{2}\geq n_{1}\geq M_{\varepsilon ,\sigma }
\end{equation*}%
for each $\sigma \in \mathfrak{S}_{\mathbb{N}^{\ast }}$. Consequently, the
series $\sum_{n}b_{n}^{\ast }\left( x_{n}\right) a_{n}$ is unconditionally
convergent in $E$. Hence $z_{\infty }:=\left( b_{n}^{\ast }\left(
x_{n}\right) a_{n}\right) _{n\in \mathbb{N}^{\ast }}$ belongs to $Z_{E}$.
Then 
\begin{equation*}
\Vert z_{r}-z_{\infty }\Vert _{W_{E}}\leq \frac{\varepsilon }{4},\; \;r\geq
N_{\varepsilon }
\end{equation*}%
Consequently $\left( W_{E},\Vert .\Vert _{W_{E}}\right) $ is a Banach space.

$\square $

\textbf{Claim 3 }: \textit{The sequence }$\left( A_{n}\right) _{n\in 
\mathbb{N}
^{\ast }}:=\left( \left( \delta _{nr}a_{n}\right) _{r\in \mathbb{N}^{\ast
}}\right) _{n\in 
\mathbb{N}
^{\ast }}$\textit{is a \ sequence of elements of }$W_{E}$\textit{\ and is an
unconditional Schauder basis of }$\left( W_{E},\Vert .\Vert _{W_{E}}\right) $%
\textit{.}

\textbf{Proof of the claim}\newline
Since $b_{n}^{\ast }\neq 0_{E^{\ast }}$ for each $n\in \mathbb{N}^{\ast }$,
there exists, for each $n\in \mathbb{N}^{\ast },$ $y_{n}\in E$ such that $%
b_{n}^{\ast }(y_{n})=1.$ It is clear that $A_{n}:=\left( b_{r}^{\ast
}(\delta _{n,r}y_{n})a_{r}\right) _{r\in \mathbb{N}^{\ast }}\in W_{E}$ for
each $n\in \mathbb{N}^{\ast }.$ Let $z:=\left( b_{n}^{\ast
}(x_{n})a_{n}\right) _{n\in \mathbb{N}^{\ast }}\in W_{E}$. We have for every 
$n\in \mathbb{N}^{\ast }$ and $\sigma _{0}\in \mathfrak{S}_{\mathbb{N}^{\ast
}}:$ 
\begin{equation*}
\left \Vert z-\overset{n}{\underset{r=1}{\sum }}b_{\sigma _{0}(r)}^{\ast
}(x_{\sigma _{0}(r)})A_{\sigma _{0}(r)}\right \Vert _{Z_{E}}=\underset{%
\sigma \in \mathfrak{S}_{\mathbb{N}^{\ast }},\text{ }r\in 
\mathbb{N}
^{\ast }}{\sup }\left \Vert \underset{1\leq j\leq r\text{ , }\sigma (j)\in 
\mathbb{N}^{\ast }\backslash \left \{ \sigma _{0}(1),...,\sigma
_{0}(n)\right \} }{\sum }b_{\sigma (j)}^{\ast }(x_{\sigma (j)})a_{\sigma
(j)}\right \Vert _{E}
\end{equation*}%
Let us recall that there exists for each $\varepsilon >0$, an integer $%
k(\varepsilon )\in \mathbb{N}^{\ast }$ such that: 
\begin{equation*}
\left \Vert \underset{j\in L}{\sum }b_{j}^{\ast }(x_{j})a_{j}\right \Vert
_{E}\leq \varepsilon
\end{equation*}%
for each $L\in \mathcal{D}_{\mathbb{N}^{\ast }}$ such that $\min (L)\geq
k(\varepsilon )+1$. It follows that we have for each integer $n\in \mathbb{N}%
^{\ast }$ such that $n\geq \underset{1\leq j\leq k(\varepsilon )}{\max }%
\left( \sigma _{0}^{-1}(j)\right) +1$: 
\begin{equation*}
\underset{\sigma \in \mathfrak{S}_{\mathbb{N}^{\ast }},\text{ }r\in \mathbb{N%
}^{\ast }}{\sup }\left \Vert \underset{1\leq j\leq r,\text{ }\sigma (j)\in 
\mathbb{N}^{\ast }\backslash \left \{ \sigma _{0}(1),...,\sigma
_{0}(n)\right \} }{\sum }b_{\sigma (j)}^{\ast }(x_{\sigma (j)})A_{\sigma
(j)}\right \Vert _{E}\leq \varepsilon
\end{equation*}%
Consequently: 
\begin{equation*}
\underset{n\rightarrow +\infty }{\lim }\left \Vert z-\overset{n}{\underset{%
r=1}{\sum }}b_{\sigma _{0}(r)}^{\ast }(x_{\sigma _{0}(r)})B_{\sigma
_{0}(r)}\right \Vert _{W_{E}}=0
\end{equation*}%
for every $\sigma _{0}\in \mathfrak{S}_{\mathbb{N}^{\ast }}$. Hence the
series $\sum b_{n}^{\ast }(x_{n})B_{n}$ is unconditionally convergent in $%
W_{E}$ to $z$. On the other hand, we have for each $n\in \mathbb{N}^{\ast }$
and $(\alpha _{1},...\alpha _{n+1})\in \mathbb{K}^{n+1}$ : 
\begin{equation*}
\alpha _{1}A_{1}+...+\alpha _{n}A_{n}=\left( b_{1}^{\ast
}(z_{r})a_{r}\right) _{r\in 
\mathbb{N}
^{\ast }}
\end{equation*}%
where :%
\begin{equation*}
z_{r}:=\left \{ 
\begin{array}{c}
\alpha _{r}y_{r}\text{ if }1\leq r\leq n \\ 
0\text{ if }r\geq n+1%
\end{array}%
\right.
\end{equation*}%
It follows that: 
\begin{align*}
& \left \Vert \alpha _{1}A_{1}+...+\alpha _{n}A_{n}\right \Vert _{W_{E}}=%
\underset{\sigma \in \mathfrak{S}_{\mathbb{N}^{\ast }},\text{ }r\in \mathbb{N%
}^{\ast }}{\sup }\left \Vert \overset{r}{\underset{j=1}{\sum }}b_{\sigma
(j)}^{\ast }(z_{\sigma (j)})a_{\sigma (j)}\right \Vert _{E} \\
& =\underset{\sigma \in \mathfrak{S}_{\mathbb{N}^{\ast }},\text{ }r\in 
\mathbb{N}^{\ast }}{\sup }\left \Vert \underset{1\leq j\leq r,\text{ }1\leq
\sigma (j)\leq n}{\sum }b_{\sigma (j)}^{\ast }(z_{\sigma (j)})a_{\sigma
(j)}\right \Vert _{E} \\
& =\underset{\sigma \in \mathfrak{S}_{\mathbb{N}^{\ast }},\text{ }r\in 
\mathbb{N}^{\ast }}{\sup }\left \Vert \underset{1\leq \sigma ^{-1}(j)\leq
r,1\leq j\leq n}{\sum }b_{j}^{\ast }(z_{j})a_{j}\right \Vert _{E} \\
& =\underset{\sigma \in \mathfrak{S}_{\mathbb{N}^{\ast }},\text{ }r\in 
\mathbb{N}^{\ast }}{\sup }\left \Vert \underset{j\in \left \{ \sigma
(1),...,\sigma (r)\right \} \cap \left \{ 1,...,n\right \} }{\sum }%
b_{j}^{\ast }(z_{j})a_{j}\right \Vert _{E} \\
& =\underset{L\in \mathcal{P}_{\left \{ 1,...,n\right \} }}{\sup }\left
\Vert \underset{j\in L}{\sum }b_{j}^{\ast }(z_{j})a_{j}\right \Vert _{E} \\
& \leq \underset{L\in \mathcal{P}_{\left \{ 1,...,n+1\right \} }}{\sup }%
\left \Vert \underset{j\in L}{\sum }b_{j}^{\ast }(z_{j})a_{j}\right \Vert
_{E} \\
& \leq \left \Vert \alpha _{1}A_{1}+...+\alpha _{n+1}A_{n+1}\right \Vert
_{W_{E}}
\end{align*}%
It follows that: 
\begin{equation*}
\left \Vert \alpha _{1}A_{1}+...+\alpha _{n}A_{n}\right \Vert _{W_{E}}\leq
\left \Vert \alpha _{1}A_{1}+...+\alpha _{n+1}A_{n+1}\right \Vert _{W_{E}}
\end{equation*}%
Consequently $\left( A_{n}\right) _{n\in \mathbb{N}^{\ast }}$ is an
unconconditionally monotone Schauder basis of $\left( W_{E},\Vert .\Vert
_{W_{E}}\right) $.

$\square $

\textbf{Claim 4 }:

1. \textit{The mapping}%
\begin{equation*}
\begin{array}{cccc}
T_{2}: & E & \rightarrow & W_{E} \\ 
& x & \mapsto & \left( b_{n}^{\ast }\left( x\right) a_{n}\right) _{n\in 
\mathbb{N}^{\ast }}%
\end{array}%
\end{equation*}%
\textit{\ is well-defined with a closed image }$T_{2}\left( E\right) $%
\textit{. Furthermore }$T_{2}$\textit{\ is an isomorphism from }$E$\textit{\
onto }$T_{2}(E)$.

2. $T_{2}(E)$\textit{\ is a complemented subspace of }$W_{E}$\textit{.}

\textbf{Proof of the claim }:

1. Since $\left( \left( a_{n},b_{n}^{\ast }\right) \right) _{n\in \mathbb{N}%
^{\ast }}$ is a besselian Schauder frame of $E$ and that $E$ is weakly
sequentially complete, it follows thanks to the proposition , that the
series $\sum b_{n}^{\ast }(x)a_{n}$ is unconditionally convergent in $E$,
for each $x\in E$. It follows that the sequence $\left( \left( b_{n}^{\ast
}(x)a_{n}\right) \right) _{n\in \mathbb{N}^{\ast }}$ belongs to $W_{E}$ for
each $x\in E$. Consequently the mapping $T_{2}$ is well defined. It is clear
that $T_{2}$ is a linear operator. Furthermore, since $\left( \left(
a_{n},b_{n}^{\ast }\right) \right) _{n\in \mathbb{N}^{\ast }}$ is a
besselian Schauder frame of $E$, we have for all $x\in E$ : 
\begin{align*}
& \left \Vert T_{2}(x)\right \Vert _{Z_{E}}=\underset{\sigma \in \mathfrak{S}%
_{\mathbb{N}^{\ast }},\text{ }n\in \mathbb{N}^{\ast }}{\sup }\left \Vert 
\overset{n}{\underset{j=1}{\sum }}b_{\sigma (j)}^{\ast }(x)a_{\sigma
(j)}\right \Vert _{E} \\
& =\underset{\sigma \in \mathfrak{S}_{\mathbb{N}^{\ast }},\text{ }n\in 
\mathbb{N}^{\ast }}{\sup }\left( \underset{y^{\ast }\in \mathbb{B}_{E^{\ast
}}}{\sup }\left \vert \overset{n}{\underset{j=1}{\sum }}b_{\sigma (j)}^{\ast
}(x)y^{\ast }(a_{\sigma (j)})\right \vert \right) \\
& \leq \underset{\sigma \in \mathfrak{S}_{\mathbb{N}^{\ast }}}{\sup }\left( 
\underset{y^{\ast }\in \mathbb{B}_{E^{\ast }}}{\sup }\overset{+\infty }{%
\underset{n=1}{\sum }}\left \vert b_{\sigma (n)}^{\ast }(x)y^{\ast
}(a_{\sigma (n)})\right \vert \right) \\
& \leq \underset{y^{\ast }\in \mathbb{B}_{E^{\ast }}}{\sup }\overset{+\infty 
}{\underset{n=1}{\sum }}\left \vert b_{n}^{\ast }(x)y^{\ast }(a_{n})\right
\vert \\
& \leq \mathcal{L_{\mathcal{F}}}\Vert x\Vert _{E}
\end{align*}%
Consequently, $T_{2}$ is a bounded linear operator. We have also for each $%
x\in E$: 
\begin{align*}
\Vert T_{2}(x)\Vert _{W_{E}}& \geq \underset{n\in \mathbb{N}^{\ast }}{\sup }%
\left( \left \Vert \overset{n}{\underset{j=1}{\sum }}b_{j}^{\ast
}(x)a_{j}\right \Vert _{E}\right) \\
& \geq \Vert x\Vert _{E}
\end{align*}%
Consequently we have :%
\begin{equation}
\left \Vert x\right \Vert _{E}\leq \left \Vert T_{2}(x)\right \Vert
_{Z_{E}}\leq \mathcal{L}_{\mathcal{F}}\left \Vert x\right \Vert _{E},\text{ }%
x\in E  \label{f}
\end{equation}%
It follows, from (\ref{f}), that $T_{2}\left( E\right) $ is a closed
subspace of $E$ and that $T_{2}$ is an isomorphism from $E$ onto $%
T_{2}\left( E\right) $.

2. Let us consider the mapping $\widetilde{\varrho }:W_{E}\longrightarrow
W_{E}$ defined for each $z:=\left( b_{n}^{\ast }(x_{n})a_{n}\right) _{n\in 
\mathbb{N}^{\ast }}\in W_{E}$ by the formula :%
\begin{equation*}
\widetilde{\varrho }(z):=\left( b_{n}^{\ast }\left( \underset{m=1}{\overset{%
+\infty }{\sum }}b_{m}^{\ast }(x_{m})a_{m}\right) a_{n}\right) _{n\in 
\mathbb{N}^{\ast }}
\end{equation*}%
It is easy to check that $\varrho $ is a well-defined linear mapping and
that : 
\begin{equation*}
\left \{ 
\begin{array}{c}
\widetilde{\varrho }(z)=T_{2}\left( \left( \underset{m=1}{\overset{+\infty }{%
\sum }}b_{m}^{\ast }(x_{m})a_{m}\right) \right) \in T_{2}\left( E\right) ,%
\text{ }z\in W_{E} \\ 
\varrho (z)=z\text{, }z\in T_{2}(E)%
\end{array}%
\right.
\end{equation*}%
It follows that $\widetilde{\varrho }(W_{E})=T_{2}(E).$

We have also :\newline
\begin{eqnarray*}
&&\left \Vert \widetilde{\varrho }(z)\right \Vert _{W_{E}}=\left \Vert
T_{2}\left( \left( \underset{m=1}{\overset{+\infty }{\sum }}b_{m}^{\ast
}(x_{m})a_{m}\right) \right) \right \Vert _{Z_{E}} \\
&\leq &\mathcal{L}_{\mathcal{F}}\left \Vert \underset{m=1}{\overset{+\infty }%
{\sum }}b_{m}^{\ast }(x_{m})a_{m}\right \Vert _{E} \\
&\leq &\mathcal{L}_{\mathcal{F}}\underset{n\in 
\mathbb{N}
^{\ast }}{\sup }\left \Vert \underset{j=1}{\overset{n}{\sum }}b_{j}^{\ast
}(x_{j})a_{j}\right \Vert _{E} \\
&\leq &\mathcal{L}_{\mathcal{F}}\left \Vert z\right \Vert _{Z_{E}}
\end{eqnarray*}

It follows that $\widetilde{\varrho }:W_{E}\rightarrow W_{E}$ is a bounded
linear operator.

By a similar way as for the proof of the direct part of the theorem $5.1.$,
we show that $\widetilde{\varrho }\circ \widetilde{\varrho }=\widetilde{%
\varrho }.$

Consequently $\widetilde{\varrho }$ is a bounded projection from $W_{E}$
onto $T_{2}(E)$. It follows that the closed subspace $T_{2}(E)$ of $W_{E}$,
is a complemented subspace of $W_{E}$.

$\square $

\textbf{End of the proof of the theorem :}

$E$ is isomorphic to the closed complemented subspace $T_{2}(E)$ of $W_{E}$
which has a monotone unconditional Schauder basis $\left( A_{n}\right)
_{n\in 
\mathbb{N}
^{\ast }}$. But there exists, according to (\cite{pel1}, page 211,
proposition 1), a Banach space $W$\ and an isomorphism $T_{3}:W_{E}%
\rightarrow W$ such that $E$ is isometrically isomorphic to the subspace $%
T_{3}\left( T_{2}\left( E\right) \right) .$ It follows that $\left(
T_{3}\left( A_{n}\right) \right) _{n\in 
\mathbb{N}
^{\ast }}$ is an unconditional Schauder basis of the Banach space $W$ and
that $T_{3}\left( T_{2}\left( E\right) \right) $ is a complemented subspace
of $W.$ Hence we achieve the proof of the theorem.
\end{proof}

\textbf{Corollary 5.1. }\textit{The space }$L_{1}\left( \left[ 0,1\right]
\right) $\textit{\ has no besselian Schauder frame.}

\begin{proof}
For the sake of contradiction, we assume that $L_{1}([0,1])$ has a besselian
Schauder frame. It is well known that the Banach space $L_{1}([0,1])$ is
weakly sequentially complete ((\cite{lin02}, page 31), (\cite{dun}, IV.8.6,
pages 290, 291, theorem 6)). It follows, thanks to the theorem $5.2.,$ that $%
L_{1}([0,1])$ is isomorphic to a subspace of a Banach space with an
unconditional basis. But this consequence contradict a well-known result (%
\cite[page:24]{lin01}, page $24$, proposition 1.d.1.) on the space $%
L_{1}([0,1])$. Consequently the Banach space $L_{1}([0,1])$ has no besselian
Schauder frame.
\end{proof}

\textbf{Corollary 5. 2.}

1. \textit{There exists a separable Banach space }$\mathcal{B}$\textit{\
with a Schauder basis\ such that a Banach space }$E$\textit{\ has a Schauder
frame if and only if }$E$\textit{\ is isomorphic to a complemented subspace
\ of }$\mathcal{B}.$\textit{\ }

2. \textit{There exists a separable Banach space} $\mathcal{U}$\textit{\
with an unconditional basis such that each weakly sequentially complete
Banach space }$E$\textit{\ which has a besselian Schauder frame is
isomorphic to a complemented subspace of $\mathcal{U}$.}

\begin{proof}
1. a. Assume that $E$ has a Schauder frame. Then there exist thanks to the
theorem 5.1., a Banach space $Z$ \ which has a Schauder basis and a bounded
linear operator $\varphi :E\longrightarrow Z$ such that $\varphi $ is an
isomorphism from $E$ onto $\varphi (E)$ and $\varphi (E)$ is a closed
complemented subspace in $Z$ . But we know, thanks to ((\cite{pel1}), page
248, Corollary 1.), that there exists a universal separable Banach space $%
\mathcal{B}$ with a Schauder basis such that every Banach space with a
Schauder basis is isomorphic to a closed subspace of $\mathcal{B}$ which is
complemented in $\mathcal{B}.$ Hence there exists a linear operator $\Phi :$ 
$Z\rightarrow \mathcal{B}$ such that $\Phi $ is an isomorphism from $Z$ onto 
$\Phi (Z)$ and $\Phi (Z)$ is both a closed subspace of $\mathcal{B}$ and a
complemented subspace in $\mathcal{B}.$\ So there exists a closed subspace $%
V_{0}$ of $Z$ and a closed subspace $V_{1}$ of $\mathcal{B}$ such that :%
\begin{equation*}
\left \{ 
\begin{array}{c}
\varphi (E)\oplus V_{0}=Z\text{ } \\ 
\Phi (Z)\oplus V_{1}=\mathcal{B}\newline
\end{array}%
\right.
\end{equation*}

It follows, from the assumptions on $\varphi $ and $\Phi $, that $\left(
\Phi \circ \varphi \right) (E)$ is closed in $\mathcal{B}$, $\Phi (V_{0})$
is closed in $\mathcal{B}\ $and 
\begin{equation*}
\Phi (Z)=\left( \Phi \circ \varphi \right) (E)\oplus \Phi (V_{0})
\end{equation*}%
Hence we have : 
\begin{equation*}
\mathcal{B}=\left( \Phi \circ \varphi \right) (E)\oplus \Phi (V_{0})\oplus
V_{1}
\end{equation*}%
It follows that $E$ is isomorphic to $\left( \Phi \circ \varphi \right) (E)$
which is a closed subspace of $\mathcal{B}$ and complemented in $\mathcal{B}$%
.

b. It is clear that if $E$ is isomorphic to a complemented subspace of $%
\mathcal{B}$ then, thanks to the theorem 5.2., $E$ will have a Schauder
frame.

2. Assume that $E$ is weakly sequentially complete and that $E$ has a \
besselian Schauder frame. Then there exist a Banach space $W$\ wich has an
unconditional Schauder basis $E$\ , a bounded linear operator $\psi
:E\longrightarrow W$ such that $\psi $ is an isomorphism from $E$ onto $\psi
(E)$ and $\psi (E)$ is a closed complemented subspace in $W$. But we know,
thanks to (\cite{pel2}, page 248, Corollary 1.), that there exists a
universal separable Banach space $\mathcal{U}$ with an unconditional
Schauder basis such that every Banach space with an unconditional Schauder
basis is isomorphic to a closed subspace of $\mathcal{U}$ which is
complemented in $\mathcal{U}.$ Hence there exists a linear operator $\Psi :$ 
$W\rightarrow \mathcal{U}$ such that $\Psi $ is an isomorphism from $W$ onto 
$\Psi (W)$ and $\Psi (W)$ is both a closed subspace of $\mathcal{U}$ and a
complemented subspace in $\mathcal{U}.$\ So there exist a closed subspace $%
V_{2}$ of $W$ and a closed subspace $V_{3}$ of $\mathcal{U}$\ such that :%
\begin{equation*}
\left \{ 
\begin{array}{c}
\psi (E)\oplus V_{2}=W \\ 
\Psi (Z)\oplus V_{3}=\mathcal{U}%
\end{array}%
\right.
\end{equation*}

It follows, from the assumptions on $\psi $ and $\Psi $, that $\left( \Psi
\circ \psi \right) (E)$ is closed in $\mathcal{U}$, $\Psi (V_{2})$ is closed
in $\mathcal{U}\ $and 
\begin{equation*}
\Psi (Z)=\left( \Psi \circ \psi \right) (E)\oplus \Psi (V_{2})
\end{equation*}%
Hence we have : 
\begin{equation*}
\mathcal{U}=\left( \Psi \circ \psi \right) (E)\oplus \Psi (V_{2})\oplus V_{3}
\end{equation*}%
It follows that $E$ is isomorphic to $\left( \Psi \circ \psi \right) (E)$
which is a closed subspace of $\mathcal{U}$ and complemented in $\mathcal{U}%
. $
\end{proof}

\bigskip

\bigskip

\bigskip

\bigskip


\begin{thebibliography}{99}
\bibitem{ald} A. Aldroubi, Q. Sun and W. Tang, \emph{p-frames and shift
invariant subspaces of L}$^{\emph{p}},$\emph{\ }J. Fourier Anal. Appl. 7(1) $%
(2001),$ 1-22.

\bibitem{ari} M. L. Arias, G. Corach, M. Pacheco, \emph{Characterisation of
Bessel sequences, \ }Extr. Math. 22 1 $(2015)$ 1-13.

\bibitem{bak} D. Bakic, T. Beric, \emph{Finite extension of Bessel
multipliers, }J. Math, Anal. Appl. 9 $(2015)$ 1-13.

\bibitem{cas1} P. G. Cassaza, \emph{The art of frame \ theory, }Taiwanese J.
Math. Phys. \ 27(5) $(1986)$, 1271-1283.

\bibitem{cas2} P. G. Cassaza, O. C hristensen and D. Stoeva, \emph{Frame
expansions in separable Banach spaces, }J. Math. Anal. Appl., 307(2) $(2005)$%
, 710-723.\emph{\ }

\bibitem{cas3} P.G. Casazza and N. Leonhard, \emph{Classes of finite equal
norm Parseval frames,} Contemp. Math. 451 $(2008)$ 11--3.

\bibitem{chr1} O. Christensen and D. T. Stoeva, \emph{p-frames in separable
Banach spaces, }Adv. Comput. Math., 18(2-4) $(2003)$, 117-1 26.\emph{\ }

\bibitem{chr2} O. Christensen, H. O. Kim and R. Y. Kim, \emph{Extensions of
Bessel sequences to dual pairs of frames,} Appl. Comput. Harmon. Anal. 34 $%
(2013)$, 224--23.

\bibitem{chr3} Ole Christensen, \emph{An Introduction to Frames and Riesz
Bases}, Birkh\"{a}user Basel $(2016).$Second Edition, (Applied and Numerical
Harmonic Analysis).

\bibitem{dau} I. Daubechies, A. Grossmann and Y. Meyer, \emph{Painless
non-orthogonal expansions,} J. Math. Phys. , New York, $1980.$

\bibitem{deh} M. A. Dehgan, M. Mesbah, \emph{Operators, frames and
convergence of sequences of Bessel sequences, }U. P. B. Sci. Bull. Series A,
77 1 $(2015)$ 75-86.\emph{\ }

\bibitem{duf} Duffin R. J. and Schaeffer A. C., \  \emph{A class of
nonharmonic Fourier series}, Trans. Amer. Math. Soc. 72 $(1952),$ 341-366.

\bibitem{dun} Dunford, N., Schwartz, J. : \emph{Linear operators, Vol. I.}
New York : Interscience $1958.$

\bibitem{fei} H. G. Feichtinger and K. Gr\"{o}chenig, \emph{A unified
approach to atomic decompositionsvia integrable group representations,} In
book : Function spaces and applications, Springer Berlin Heidelberg, $%
(1988), $ 52-73.

\bibitem{gab} Gabor, D. : \emph{Theory of communication,} J. Inst. Electr.
Eng. 93, \ 429-457 $(1946).$

\bibitem{gro} K. Gr\"{o}chenig. \emph{Describing functions: atomic
decompositions versus frames}, Monatsh. Math., Soc.,112 $(1991)$, no.1,1-42.%
\emph{\ }

\bibitem{jia} R. Q. A. Jia, Bessel sequences in Sobolev spaces, Appl.
Comput.Harmonic. Anal. 20 \ $(2006)$ 298-311.

\bibitem{lar1} P. G. Casazza, D. Han and D. R. Larson, \emph{Frames for
Banach spaces,} Contemp. Math., 247 $(1999)$, 149-182.

\bibitem{lar2} D. Han, D. R. Larson, \emph{Frames, bases and group
representation,} Mem. Amer. Math. Soc. 147(697) $(2000)$, 1-91.

\bibitem{li} D.F. Li and W.C. Sun, \emph{Expansion of frames to tight frames}%
, Acta Math. Sin. (Engl. Ser.) 25 $(2009),$ 287-292.

\bibitem{lin01} J. Lindenstraus, J. Tzafriri. \emph{Classical Banach spaces
I. Sequence spaces.} Reprint of the $1977$ edition. Springer.

\bibitem{lin02} J. Lindenstraus, J. Tzafriri. \emph{Classical Banach spaces
II. Function spaces.} Reprint of the $1979$ edition. Springer.

\bibitem{lin1} J. Lindenstraus and J. Tzafriri. \emph{On the complemented
subspaces problem.} Israel J. Math. 9, 263--269 $(1971)$.

\bibitem{meg} R. E. Megginson, \emph{An Introduction to Banach Space Theory,}
Graduate Texts in Mathematics, vol. 183 (Springer, New York, $1998$)

\bibitem{pel1} Pe\l czy\'{n}ski, A. $(1960).$ \emph{Projections in certain
Banach spaces.} Studia Mathematica, 2(19), 209-228.

\bibitem{pel2} Pe\l czy\'{n}ski, A.. \emph{"Universal bases."} Studia
Mathematica 32.3 $(1969)$ : 247-268.

\bibitem{pel3} Pe\l czy\'{n}ski, A.. \emph{"Any separable Banach space with
the bounded approximation property is a complemented subspace of a Banach
space with a basic."} Studia Mathematica 40.3 $(1971)$ : 239-243.

\bibitem{rah} A. Rahimi, P. Balazs, \emph{Multipliers for p-Bessel sequences
in Banach spaces,} Integral Equations Operator Theory, 68 $(2010)$ 193-205.

\bibitem{rya} Ryan, Raymond A. $(2002).$ \emph{Introduction to Tensor
Products of Banach Spaces.} Springer Monographs in Mathematics. London New
York: Springer.

\bibitem{sto} D. T. Stoeva, \emph{Frames and Bases in Banach spaces, }PhD
thesis, $2005.$

\bibitem{you} R. Young, \emph{An introduction to Nonharmonic Fourier Series, 
}Academic Press,\emph{\ New York, }$1980.$

\bibitem{woj} P. Wojtaszczyk, \emph{Banach spaces for analysts.} (Cambridge
studies in advanced mathematics 25), Cambridge Univerisity Press $1991.$
\end{thebibliography}
\end{document}